\documentclass[11pt]{amsart}

\usepackage{amsmath,amsthm, amscd, amssymb, amsfonts}
\usepackage{epsfig}
\usepackage{verbatim}

\newcommand{\otb}{{\overline{\otimes}}}
\newcommand{\otk}{{\otimes}}
\newcommand{\Mo}{{\mathcal M}}

\newcommand{\nic}{{\mathfrak B}}

\newcommand{\Ss}{{\mathcal S}}
\newcommand{\ot}{{\otimes}}

\newcommand{\kuc}{{\mathcal K}}
\newcommand{\Ac}{{\mathcal A}}
\newcommand{\Bc}{{\mathcal B}}

\newcommand{\ca}{{\mathcal C}}

\newcommand{\Do}{{\mathcal D}}

\newcommand{\YD}{{\mathcal YD}}

\newcommand{\ku}{{\Bbbk}}

\newcommand{\Z}{{\mathbb Z}}
\newcommand{\Na}{{\mathbb N}}

\newcommand{\uno}{{\mathbf 1}}

\newcommand{\id}{\mbox{\rm id\,}}

\newcommand\stab{\operatorname{Stab}}

\newcommand\Rep{\operatorname{Rep}}

\newcommand{\End}{\operatorname{End}}

\newcommand{\ad}{\mbox{\rm ad\,}}

\newcommand{\gr}{\mbox{\rm gr\,}}

\renewcommand{\_}[1]{\mbox{$_{\left( #1 \right)}$}}

\theoremstyle{plain}

\numberwithin{equation}{section}

\newtheorem{teo}{Theorem}[section]

\newtheorem{lema}[teo]{Lemma}

\newtheorem{cor}[teo]{Corollary}

\newtheorem{prop}[teo]{Proposition}

\newtheorem{claim}{Claim}[section]

\theoremstyle{definition}

\newtheorem{defi}[teo]{Definition}

  \newtheorem{exa}[teo]{Example}

\newtheorem{exas}[teo]{Examples}

\theoremstyle{remark}

\newtheorem{rmk}[teo]{Remark}

\def\pf{\begin{proof}}

\def\epf{\end{proof}}

\theoremstyle{remark}

\def\s{\mathbb{S}}
\def\mO{\mathcal{O}}
\def\mK{\mathcal{K}}
\def\mA{\mathcal{A}}
\def\mJ{\mathcal{J}}
\def\mT{\mathcal{T}}
\def\mR{\mathcal{R}}
\def\B{\mathfrak{B}}

\def\mQ{\mathcal{Q}}
\def\mD{\mathcal{D}}
\def\mH{\mathcal{H}}
\def\mL{\mathcal{L}}

\def\hq{\mH(\mQ)}
\def\N{\mathbb{N}}

\newcommand{\ydg}{{}^{G}_{G}\mathcal{YD}}
\newcommand\sg{\operatorname{sgn}}

\def\eps{\epsilon}
\def\veps{\varepsilon}
\newcommand\can{\operatorname{can}}
\begin{document}

\title[Module categories over pointed Hopf algebras over $\s_3$ and
$\s_4$]{Representations of the category of modules over pointed Hopf algebras
over $\s_3$ and $\s_4$}
\author[Garc\'\i a Iglesias and Mombelli]{Agust\'\i n Garc\'\i a Iglesias and
Mart\'\i n Mombelli }

\address{ Facultad de Matem\'atica, Astronom\'\i a y F\'\i sica,
Universidad Nacional de C\'ordoba, CIEM, Medina Allende s/n,
(5000) Ciudad Universitaria, C\'ordoba, Argentina}
\email{aigarcia@famaf.unc.edu.ar,\newline \indent \qquad\qquad\qquad\quad
martin10090@gmail.com
\newline \indent\emph{URL:}\/http://www.mate.uncor.edu/$\sim$aigarcia,\newline
\indent \quad\quad  http://www.mate.uncor.edu/$\sim$mombelli}

\thanks{ {\sc keywords}: Tensor categories, module categories, pointed Hopf
algebras.\newline
 The work was partially supported by CONICET,
FONCyT-ANPCyT, Secyt (UNC), Mincyt (C\'ordoba)}

\begin{abstract}
We  classify exact indecomposable  module categories over the representation
category of all non-trivial Hopf algebras with coradical $\s_3$ and  $\s_4$. As
a byproduct, we compute all its Hopf-Galois extensions  and we show that these
Hopf algebras are cocycle deformations of their graded versions.
\end{abstract}

\date{\today}
\maketitle

\hspace{3cm}

\section{Introduction}

Given a tensor category $\ca$, an  \emph{exact module category} \cite{eo} over
$\ca$ is an Abelian category $\Mo$ equipped with a biexact functor
$\otimes:\ca\times \Mo\to \Mo$ subject to natural associativity and unity
axioms, such that for any projective object $P\in \ca$ and any $M\in \Mo$ the
object $P\otimes M$ is again projective.

\medbreak

Exact module categories, or \emph{representations} of $\ca$, are  very
interesting objects to consider. They are implicitly present in many areas of
mathematics and mathematical physics such as  subfactor theory \cite{BEK},
affine Hecke algebras \cite{BO}, extensions of vertex algebras
\cite{KO}, \cite{HuKo}, Calabi-Yau algebras \cite{Gi} and conformal field
theory, see for example \cite{BFRS}, \cite{FS}, \cite{CS1}, \cite{CS2}.

Module categories have been  used in the study of fusion categories \cite{ENO1},
\cite{ENO1}, and in the theory of (weak) Hopf algebras \cite{O1}, \cite{M1},
\cite{N}.

\medbreak

The classification of exact module categories over a fixed finite tensor
category $\ca $ has been undertaken by several authors:
\begin{itemize}
  \item[1.] When $\ca$ is the semisimple quotient of $U_q(\mathfrak{sl}_2)$
\cite{KO}, \cite{eo2},
        \item[2.] over the tensor categories of representations of finite
supergroups \cite{eo},
  \item[3.] over $\Rep(D(G))$, $D(G)$ the Drinfeld double of a finite group $G$
\cite{O2},
  \item[4.] over the tensor category of representations of the Lusztig's small
quantum group $u_q(\mathfrak{sl}_2)$ \cite{M1},
  \item[5.] and more generally over $\Rep(H)$, where $H$ is a lifting of a
quantum linear space \cite{M2}.
\end{itemize}

\smallbreak

The main goal of this paper is the classification of exact module categories
over the representation category of any non-trivial (i.e. different from the
group algebra) finite-dimensional Hopf algebra with coradical $\ku \s_3$ or
$\ku \s_4$.

\medbreak

Finite-dimensional Hopf algebras with coradical $\ku \s_3$ or  $\ku \s_4$ were
classified in \cite{AHS}, \cite{GG}, respectively. For all these Hopf algebras the associated graded Hopf algebras $\gr H$ are of the form $\nic(X,q)\#
\ku \s_n$, $n=3,4$ where $X$ is a finite set equipped with a map $\rhd:X\times
X\to X$ satisfying certain axioms that makes it into a \textit{rack} and
$q:X\times X\to \ku^{\times}$ is a 2-cocycle. We have the following result:

Let $n=3,4$ and let $\Mo$ be an exact indecomposable module category over
$\Rep(\nic(X,q)\# \ku \s_n)$, then there exist
\begin{itemize}
  \item a subgroup $F< \s_n$ and a 2-cocycle $\psi\in Z^2(F,\ku^{\times})$,
  \item a subset $Y\subseteq X$ invariant under the action of $F$,
  \item a family of scalars $\{\xi_C\}$ compatible (Definition \ref{def:compatible}) with $(F, \psi, Y)$,
\end{itemize}
such that $\Mo\simeq {}_{\Bc(Y, F, \psi, \xi)}\Mo$, where $\Bc(Y, F, \psi, \xi)$
is a left $\nic(X,q)\# \ku \s_n$-comodule algebra constructed from data $(Y, F,
\psi, \xi)$. We also show that if $H$ is a finite-dimensional Hopf algebra with
coradical $\ku \s_3$ or  $\ku \s_4$ then $H$ and $\gr H$ are cocycle
deformations of each other. This implies that there is a bijective
correspondence between module categories over $\Rep(H)$ and $\Rep(\gr H)$.

\medbreak

The content of the paper is as follows. In Section \ref{modc} we recall the
basic results on module categories over finite-dimensional Hopf algebras. We
recall the main result of \cite{M2} that gives an isomorphism between
Loewy-graded comodule algebras and a semidirect product of a twisted group
algebra and a homogeneous coideal subalgebra inside the Nichols algebra. In
Section \ref{equivariant:classes} we show how to distinguish Morita
equivariant classes of comodule algebras over pointed Hopf algebras.

\medbreak

In Section \ref{pointed} we recall the definition of a rack $X$ and a ql-datum
$\mQ$, and how to construct (quadratic approximations to) Nichols algebras
$\widehat{\B_2}(X,q)$ and pointed Hopf algebras $\hq$ out of them. In
particular, we recall a presentation of all finite-dimensional Hopf algebras
with coradical $\ku\s_3$, $\ku\s_4$. In Section \ref{sec:coideal}, we give a
classification of connected homogeneous left coideal subalgebras of
$\widehat{\B_2}(X,q)$ and also a presentation by generators an relations.

\medbreak

In Section \ref{main:section} we introduce a family of comodule algebras large
enough to classify module categories. We give an explicit Hopf-biGalois
extension over $\widehat{\B_2}(X,q)\# \ku \s_n $, $n\in\N$, and a lifting $\hq$,
proving that there is a bijective correspondence between module categories over
$\Rep(\widehat{\B_2}(X,q)\# \ku \s_n)$ and $\Rep(\hq)$, $n=3,4$. In particular we obtain that any pointed Hopf
algebra over $\s_3$ or $\s_4$ is a cocycle deformation of its associated
graded algebra, a result analogous to a theorem of Masuoka for abelian groups \cite{Ma}.
Finally, the classification of module categories over
$\Rep(\widehat{\B_2}(X,q)\# \ku \s_n)$ is presented in this section and as a
consequence all Hopf-Galois objects over $\widehat{\B_2}(X,q)\# \ku \s_n $ are
described.

\subsection*{Acknowledgments.} We thank N. Andruskiewitsch for suggesting us this project and for his comments  that  improved the presentation of the paper.

\section{Preliminaries and notation}

We shall denote by $\ku$ an algebraically closed field of
characteristic zero. The tensor product over the field $\ku$ will be denoted by
$\otimes$. All vector spaces, algebras and categories
will be considered over $\ku$. For any algebra $A$,  ${}_A\Mo$
will denote the category of finite-dimensional left $A$-modules.

\medbreak

The symmetric group on $n$ letters is denoted by $\s_n$ and  by $\mO_j^n$ we
shall denote the conjugacy class of all $j$-cycles in
$\s_n$. For any group $G$, a 2-cocycle $\psi\in Z^2(G,\ku^{\times})$ and any
$h\in G$ we shall denote $ \psi^h(x,y)=\psi(h^{-1}xh, h^{-1}yh)$ for all $x,y\in
G$.

\medbreak

If $H$ is a Hopf algebra, a 2-cocycle $\sigma$ in $H$ is a convolution invertible linear map $\sigma:H\times H\to\ku$ such that
\begin{equation}\label{eqn:2cociclo}
\sigma(x_{(1)}, y_{(1)}) \sigma(x_{(2)}y_{(2)},z) =
\sigma(y_{(1)},z_{(1)}) \sigma(x,y_{(2)} z_{(2)})
\end{equation}
and $\sigma(x,1) = \sigma(1,x) = \varepsilon(x)$, for every $x,y,z \in H$.
The set of 2-cocycles in $H$ is denoted by $Z^2(H)$.

\medbreak

If $A$ is an $H$-comodule algebra via $\lambda:A\to H\otk A$, we
shall say that a (right) ideal $J$ is $H$-costable if
$\lambda(J)\subseteq H\otk J$. We shall say that $A$ is (right)
$H$-simple if there is no nontrivial (right) ideal $H$-costable in
$A$.

\medbreak

If $H=\bigoplus H(i)$ is a coradically graded Hopf algebra we shall
say that a left coideal subalgebra $K\subseteq H$ is
\emph{homogeneous} if $K=\bigoplus K(i)$ is graded as an algebra
and, for any $n$, $K(n)\subseteq H(n)$ and
$\Delta(K(n))\subseteq \bigoplus_{i=0}^n \, H(i)\otk K(n-i)$. $K$ is said to be
\emph{connected} if $\mK\cap H(0)=\ku$.

If $H=\nic(V)\# \ku G$,
where $V$ is a Yetter-Drinfeld module over $G$ and $K\subseteq H$
is a coideal subalgebra, we shall denote by $\stab K$ the biggest
subgroup of $G$ such that the adjoint action of $\stab K$ leaves $K$ invariant.

\medbreak

If $H$ is a finite-dimensional Hopf algebra then $H_0\subseteq H_1
\subseteq \dots \subseteq H_m=H$ will denote the coradical
filtration. When $H_0\subseteq H$ is a Hopf subalgebra then the
associated graded algebra $\gr H$ is a coradically graded Hopf
algebra. If $(A, \lambda)$ is a left $H$-comodule algebra, the
coradical filtration on $H$ induces a filtration on $A$, given by
$A_n=\lambda^{-1}(H_n\otk A)$. This filtration is called the
\emph{Loewy series }on $A$.

The associated graded algebra $\gr A$ is a left $\gr H$-comodule algebra. The
algebra $A$ is right $H$-simple if and only if $\gr A$ is right $\gr H$-simple,
see \cite[Section 4]{M1}.

\section{Representations of tensor categories}\label{modc}

Given $\ca=(\ca, \ot, a, \uno)$ a tensor category, a \emph{module category} over
$\ca$ or a \emph{representation} of $\ca$ is an
Abelian category $\Mo$ equipped with an exact bifunctor $\otb:
\ca \times \Mo \to \Mo$ and natural associativity and unit
isomorphisms $m_{X,Y,M}: (X\otimes Y)\otimes M \to X\otimes
(Y\otimes M)$, $\ell_M: \uno \otimes M\to M$ satisfying natural
associativity and unit axioms, see \cite{eo}, \cite{O1}. We shall
assume, as in \cite{eo}, that all module categories have only
finitely many isomorphism classes of simple objects.

\medbreak

A module category is {\em indecomposable} if it is not equivalent
to a direct sum of two non trivial module categories. A module
category $\Mo$ over a finite tensor category $\ca$ is \emph{exact}
\cite{eo} if  for any projective $P\in \ca$ and any $M\in \Mo$,
the object $P\ot M$ is again projective in $\Mo$.

\medbreak

If $\Mo$ is an exact module category over $\ca$ then the dual
category $\ca^*_{\Mo}$, see \cite{eo}, is a finite tensor category. There is a
bijective correspondence between the set of equivalence classes of exact module
categories over $\ca$ and over $\ca^*_{\Mo}$, see \cite[Theorem 3.33]{eo}. This
implies that for any finite-dimensional Hopf algebra there is a bijective
correspondence between the set of equivalence classes of exact module categories
over $\Rep(H)$ and $\Rep(H^*)$.

\subsection{Module categories over pointed Hopf algebras}

We are interested in exact indecomposable module categories over the
representation category of finite-dimensional Hopf algebras. If $H$ is a Hopf
algebra and $\lambda:\Ac\to H\otk \Ac$ is a left $H$-comodule algebra the category ${}^H\Mo_{\Ac}$ is the category of finite-dimensional right $\Ac$-modules left $H$-comodules where the comodule structure is a $\Ac$-module morphism. If $\Ac'$ is another left $H$-comodule algebra the category ${}_{\Ac}^H\Mo_{\Ac'}$ is defined analogously.

The
category of finite-dimensional left $\Ac$-modules  ${}_\Ac\Mo$ is a representation
of $\Rep(H)$. The action $\otb:\Rep(H)\times {}_\Ac\Mo\to {}_\Ac\Mo$ is given by
$V\otb M=V\otk M$ for all $V\in \Rep(H)$, $M\in {}_\Ac\Mo$. The left $\Ac$-module
structure on $V\otk M$ is given by the coaction $\lambda$.

\medbreak

If $\Mo$ is an exact indecomposable module over $\Rep(H)$ then
there exists a left $H$-comodule algebra $\Ac$ right $H$-simple with trivial
coinvariants such that $\Mo\simeq {}_\Ac\Mo$ as modules over $\Rep(H)$ see
\cite[Theorem 3.3]{AM}.

\medbreak

If $\Ac$, $\Ac'$ are two right $H$-simple left $H$-comodule algebras such that the
categories ${}_\Ac\Mo$,  ${}_{\Ac'}\Mo$ are equivalent as representations over
$\Rep(H)$. Then there exists an equivariant Morita context $(P,Q,f,g)$, that is
$P\in {}_{\Ac'}^H\Mo_{\Ac}$, $Q\in {}_{\Ac}^H\Mo_{\Ac'}$ and $f:P\ot _{\Ac} Q\to
\Ac'$, $g:Q\ot_{\Ac'} P\to \Ac$ such that they are bimodule isomorphisms. Moreover, it holds that
$\Ac'\simeq \End_{\Ac}(P)$ as comodule algebras. The comodule structure on
$\End_{\Ac}(P)$ is given by $\lambda(T)=T\_{-1}\ot
T\_0$, where
\begin{equation}\label{h-comod}
\langle\alpha, T\_{-1}\rangle\,
T\_0(p)=\langle\alpha, T(p\_0)\_{-1}\Ss^{-1}(p\_{-1})\rangle\,
T(p\_0)\_0,\end{equation} for any $\alpha\in H^*$,
$T\in\End_{\Ac}(P)$, $p\in P$. See \cite{AM} for more details.

By the previous paragraph, we can see that the categories ${}^H\Mo_{\Ac}$ play a
central role in the theory. The following theorem will be of great use in the
next section.
\begin{teo}\label{freeness}
Let $H$ be a Hopf algebra and $\Ac$ a left $H$-comodule algebra, both finite dimensional. \begin{itemize}
\item\cite[Theorem 3.5]{Sk} If $\Ac$ is $H$-simple and $M \in {}^H\Mo_\Ac$, then there exists
$t\in \Na$ such that $M^t$, the direct sum of $t$ copies of $M$,
is a free $\Ac$-module.
\item\cite[Theorem 4.2]{Sk} $M \in {}^H\Mo_\Ac$ is free as $A$-module if and only if there exists a maximal ideal $J\subset \Ac$ such that $M/M\cdot J$ is free as $\Ac/J$-module.\qed
\end{itemize}
\end{teo}

The first statement of this theorem is actually present in the proof of \cite[Theorem 3.5]{Sk}. The second one will be particularly useful when the ideal $J$ is such that $\Ac/J=\ku$, since in this case $M/M\cdot J$ is automatically free.

\begin{teo}\cite[Theorem 3.3]{M2}\label{mod-over-pointed} 
Let $G$ be a finite group and let $H$ be a finite-dimensional pointed Hopf
algebra with coradical $\ku G$. Assume there is $V\in {}^{ G}_{ G}
\mathcal{YD}$ such that $\gr H=U=\nic(V)\# \ku G$. Let $\Ac$ be a left
$H$-comodule algebra  right $H$-simple with trivial coinvariants. There exists
\begin{itemize}
  \item[1.] a subgroup $F\subseteq G$,
  \item[2.] a 2-cocycle $\psi\in Z^2(F, \ku^{\times})$,
  \item[3.] a homogeneous left coideal subalgebra $\kuc=\oplus_{i=0}^m
\kuc(i)\subseteq \nic(V)$ such that $ \kuc(1)\subseteq V$ is a $\ku
G$-subcomodule invariant under the action of $F$,
\end{itemize}
such that $\gr \Ac\simeq \kuc\# \ku_{\psi} F$ as left $U$-comodule algebras.\qed
\end{teo}

The algebra structure and the left $U$-comodule structure of
$\kuc\#\, \ku_{\psi} F$ is given as follows. If $x,y\in \kuc$,
$f,g\in F$ then
\begin{align*} (x\# g) (y\# f)=x( g\cdot y)\# \psi(g,f)\, gf,\\
\lambda(x\# g)=(x\_1\# g )\ot (x\_2\# g),
\end{align*}
where the action of $F$ on $\kuc$ is the restriction of the
action of $G$ on $\nic(V)$ as an object in
${}_{G}^{G}\YD$. Observe that $F$ is necessarily a subgroup of $\stab \kuc$.

\section{Equivariant  equivalence classes of comodule
algebras}\label{equivariant:classes}

In this section we shall present how to distinguish equivalence classes of some
comodule
algebras over pointed Hopf algebras and then we will apply this result to our cases.
Much of the ideas here are already contained in \cite{M1}, \cite{M2} although
with less generality.

\medbreak

Let $\Gamma$ be a finite group and $H$ be a finite-dimensional pointed Hopf
algebra with coradical $\ku \Gamma$ and with coradical filtration $H_0\subseteq
H_1\subseteq \dots \subseteq H_m=H$. Assume there is $V\in {}^{ \Gamma}_{ \Gamma}
\mathcal{YD}$ such that $\gr H=U=\nic(V)\# \ku \Gamma$.

\medbreak

We begin with the following lemma.

\begin{lema}\label{lem:extension cociclo}
Let $\Gamma, U$ as above. Let $\sigma\in Z^2(\Gamma,\ku^\times)$ be a 2-cocycle. Then there exists a 2-cocycle $\varsigma\in Z^2(U)$ such that $\varsigma_{|\Gamma\times\Gamma}=\sigma$.
\end{lema}
\pf
Let us consider the linear map $\varsigma:U\times U\to\ku$ defined, on homogeneous elements $x,y\in U$ by
\begin{align*}
\varsigma(x,y)=\begin{cases}
               \sigma(x,y), & \text{if } x,y\in U(0);\\
0, &\text{otherwise.}
               \end{cases}
\end{align*}
Notice that $\varsigma(x,1)=\varsigma(1,x)=\eps(x)$ by definition. We have to check that, for $x\in U(m)$, $y\in U(n)$, $z\in U(k)$, $m,n,k\in\N$, \eqref{eqn:2cociclo} holds. Now, if $k>0$, the LHS of \eqref{eqn:2cociclo} is zero. Set $\Delta(z)=\sum_{i=0}^kz^i\ot z^{k-i}$, with $z^s\in U(s)$, $s=0,\dots,k$. Analogously, set $\Delta(y)=\sum_{j=0}^ny^j\ot y^{n-j}$, with $y^t\in U(t)$, $t=0,\dots,n$. Then, the RHS is
$$
\sum_{i=0}^k\sum_{j=0}^n\varsigma(x,y^{n-j}z^{k-i})\varsigma(y^j,z^i)=\varsigma(x,y^nz^k)=0,
$$
and thus \eqref{eqn:2cociclo} holds. Both sides of this equation are similarly seen to be zero if $m>0$ or $n>0$, while the case $m=n=k=0$ holds by definition of $\varsigma$. This map is convolution invertible and its inverse $\varsigma^{-1}$ is defined in an analogous manner, using $\sigma^{-1}$.
\epf

Let $\Ac, \Ac'$ be
two right $H$-simple left $H$-comodule algebras. Let $F, F'\subseteq \Gamma$ be subgroups and let $\psi \in Z^2(F,
\ku^{\times})$, $\psi' \in Z^2(F', \ku^{\times})$ be two cocycles such that
$\Ac_0=\ku_{\psi} F$ and $\Ac'_0=\ku_{\psi'} F'$.
Let $K, K'\in  \nic(V)$ be two homogeneous coideal subalgebras such that $\gr
\Ac = K\# \ku_{\psi} F$ and $\gr  \Ac' = K'\# \ku_{\psi'} F'$.

\medbreak

The main result of this section is the following.
\begin{teo}\label{Morita-equivalence} The categories ${}_\Ac\Mo$, ${}_{\Ac'}\Mo$
are equivalent as modules over $\Rep(H)$ if and only if there exists an element
$g\in \Gamma$ such that $\Ac'\simeq g\Ac g^{-1}$ as comodule algebras.
\end{teo}

\pf
Let us assume that ${}_\Ac\Mo\cong {}_{\Ac'}\Mo$
as $\Rep(H)$-modules. By \cite[Proposition 1.24]{AM} there exists an equivariant Morita context
$(P,Q,f,h)$. That is $P\in {}_{\Ac'}^H\Mo_{\Ac}$, $Q\in {}_{\Ac}^H\Mo_{\Ac'}$
and $f:P\ot _{\Ac} Q\to \Ac'$, $h:Q\ot_{\Ac'} P\to \Ac$ are bimodule
isomorphisms and $\Ac'\simeq \End_{\Ac}(P)$ as comodule algebras. The comodule
structure on $\End_{\Ac}(P)$ is given by $\lambda:\End_{\Ac}(P)\to H\otk
\End_{\Ac}(P)$, $\lambda(T)=T\_{-1}\ot
T\_0$ where
\begin{equation}\label{h-comod1} \langle\alpha, T\_{-1}\rangle\,
T\_0(p)=\langle\alpha, T(p\_0)\_{-1}\Ss^{-1}(p\_{-1})\rangle\,
T(p\_0)\_0,\end{equation} for any $\alpha\in H^*$,
$T\in\End_{\Ac}(P)$, $p\in P$.

\medbreak

For any $i=0,\dots, m$ define $P(i)=P_i/P_{i-1}$, where
$P_{-1}=0$. The graded vector space $\gr P=\oplus_{i=0}^m \, P(i)$
has an obvious structure that makes it into an object in the
category ${}^{U}\Mo_{K\# \ku_{\psi} F}$. We shall denote $\overline{\delta}:\gr
P\to U\otk \gr P$ the coaction. In particular $\gr P\in {}^{U}\Mo_K$, thus by
Theorem \ref{freeness} (ii) we have that $\gr P\simeq M\otk K$, where $M=\gr P/(\gr
P\cdot K^+)$, since $K/K^+=\ku$.

\medbreak

We have that $\overline{\delta}(\gr P\cdot K^+)\subset (U\ot \gr P)(K^+\ot 1+U\ot
K^+)$, since $K=\ku\oplus K^+$ and thus the map $\overline{\delta}$ induces a new map $\widehat{\delta}:M\to U'\ot M$, where $U'=U/UK^+U$.
Notice that $U'$ is a pointed Hopf algebra with coradical $\ku\Gamma$, since $U$ is coradically graded and the ideal $UK^+U$ is homogeneous and does not intersect $U_0$. $M$ is also a $\ku_\psi F$-module with
$\overline{m}\cdot f= \overline{m\cdot f}$, for $f\in F$, $\overline{m}\in M$.
This action is easily seen to be well defined and, moreover, $M\in {}^{U'}\Mo_{\ku_{\psi} F}$.

\medbreak

Let $\Psi\in Z^2(\Gamma,\ku^*)$ be a 2-cocycle such that $\Psi_{|F\times F}=\psi$ see \cite[Proposition III (9.5)]{Br}.
Let $\zeta\in Z^2(U')$ be such that $\zeta_{|\Gamma\times
\Gamma}=\Psi^{-1}$, as in Lemma \ref{lem:extension cociclo}.
By \cite[Lemma 2.1]{M1} there exists and equivalence of categories  ${}^{U'^{\zeta}}\Mo_{\ku F}\simeq {}^{U'}\Mo_{(\ku
F)_{\Psi}}$. By Theorem \ref{freeness} (ii) any object in ${}^{U'^{\zeta}}\Mo_{\ku F}$ is a free $\ku
F$-module. Thus there is an object $N$ in  ${}^{U/U(\ku F)^+}\Mo$ such that $\gr
P\simeq N\otk K\otk \ku_{\psi} F$. Whence $\dim P= (\dim N)(\dim \Ac)$.
Similarly we can assume that there is $s\in\Na$ such that $\dim Q= s \dim \Ac'.$

\medbreak

Using Theorem \ref{freeness} (i) there exists $t\in\Na$ such that $P^t$ is a free
right $\Ac$-module, that is, there is a  vector space $T$ such that $P^t\simeq
T\otk \Ac$, hence
\begin{align}\label{dimension1} t\, \dim N=\dim T.
\end{align}
Since $P\ot _{\Ac} Q\simeq \Ac'$ then $P^t\ot _{\Ac} Q\simeq T\otk Q\simeq
\Ac'^t$, then $s\dim T \dim \Ac'=t\dim \Ac'$ and using \eqref{dimension1} we
obtain that $s\, \dim N=1$ whence $\dim N=1$ and thus $\dim P=\dim \Ac$.

\medbreak

\begin{claim}\label{claim}
Let $n\in P_0$, then $P=n\cdot\Ac$.
\end{claim}
Notice that $P_0\neq 0$. In fact, if $P_0=0$ and $k\in \N$ is minimal with $P_k\neq 0$, then $\lambda(P_k)\subset\sum_{j=0}^k H_{k-j}\ot P_j=H_0\ot
P_k$, which is a contradiction. Let $g\in\Gamma$ be such that $\lambda(n)=g\ot n$.
Now, if $J=\{a\in\Ac\,:\, n\cdot a=0 \}$, then $J$ is a right ideal of $\Ac$.
We shall prove that $J=0$. Let $a\in J$ and write $\lambda(a)=\sum_{i=1}^n
a^i\ot a_i$, in such a way that the set $\{a^i\,:\,i=1,\dots,n \}\subset H$ is
linearly independent. Now, $\{ga^i\,:\,i=1,\dots,n \}\subset H$ is also linearly
independent and we have $0=\lambda(n\cdot a)=\sum_{i=1}^nga^i\ot n\cdot a_i$.
Thus $n\cdot a_i=0$, $\forall\,i=1,\dots, n$, that is, $\lambda(a)\in H\ot J$ and
$J$ is $H$-costable. As $\Ac$ is right $H$-simple, $J=0$. Therefore, the action
$\cdot : N\ot \Ac\to P$ is injective and since $\dim P=\dim N\dim \Ac$, the
claim follows.

\medbreak

It is not difficult to prove that the linear map $\phi:g\Ac g^{-1}\to
\End_{\Ac}(P)$ given by $ \phi(gag^{-1})(n\cdot b)=n\cdot ab$ is an isomorphism
of $H$-comodule algebras.

\medbreak

Conversely, if $\Ac'\simeq g\Ac g^{-1}$ as comodule algebras and $M\in
{}_{\Ac}\Mo$, then the set $gMg^{-1}$ has a natural structure of
$\Ac'$-module in such a way that the functor $F:{}_{\Ac}\Mo\to{}_{\Ac'}\Mo$, $M\mapsto
gMg^{-1}$ is an equivalence of $\Rep(H)$-modules.
\epf

\section{Pointed Hopf algebras over $\s_3$ and $\s_4$}\label{pointed}
In this section we describe all pointed Hopf algebras whose coradical is the
group algebra of the groups $\s_3$ and $\s_4$. These were classified in
\cite{AHS} and \cite{GG}, respectively.

\medbreak

Recall that a \emph{rack} is a pair $(X,\rhd)$,
where $X$ is a
non-empty set and $\rhd:X\times X\to X$ is a function, such that
$\phi_i=i\rhd (\cdot):X\to X$ is a bijection for all $i\in X$
satisfying:
$
i\rhd(j\rhd k)=(i\rhd j)\rhd (i\rhd k),$ for all $ i,j,k\in X.
$
See \cite{AG2} for detailed information on racks.

Let $(X,\rhd)$ be a rack. A
2-cocycle $q:X\times X\to \ku^{\times}$,
$(i,j)\mapsto q_{ij}$ is a function such that for all $\,i,j,k\in X$
$$ q_{i,j\rhd k}\, q_{j,k}=q_{i\rhd j,i\rhd
k}\, q_{i,k}.
$$

In this case it is possible to generate a braiding $c^q$ in the
vector space $\ku X$ with basis $\{x_i\}_{i\in X}$ by $c^q(x_i\otimes
x_j)=q_{ij}x_{i\rhd j}\otimes x_i,$ for all $i,j\in X$. We denote by
$\B(X,q)$ the Nichols algebra of this braided vector space.

\subsection{Quadratic approximations to Nichols algebras}\label{quadratic-app}

Let $\mJ=\oplus_{r\geq 2}\mJ^r$ be the defining ideal of the Nichols
algebra $\B(X,q)$. Next, we give a description of the space $\mJ^2$
of quadratic relations. Let $\mR$ be the set of equivalence classes
in $X\times X$ for the relation generated by $(i,j)\sim(i\rhd j,i)$.
Let $C\in\mR$, $(i,j)\in C$. Take $i_1=j$, $i_2=i$ and, recursively,
$i_{h+2}=i_{h+1}\rhd i_h$. Set $n(C)=\#C$ and $$\mR'=\Big\{C\in
\mR\,|\,\prod_{h=1}^{n(C)} q_{i_{h+1},i_h}=(-1)^{n(C)}\Big\}.$$  Let
$\mT$ be the free associative algebra in the variables
$\{T_l\}_{l\in X}$. If $C\in\mR'$, consider the quadratic polynomial
\begin{equation}\label{eqn:phiC}\phi_C = \sum_{h=1}^{n(C)}\eta_h(C)
\, T_{i_{h+1}}T_{i_h}\in \mT,\end{equation} where $\eta_1(C)=1$ and
$\eta_h(C)=(-1)^{h+1}q_{{i_2i_1}}q_{{i_3i_2}}\ldots
  q_{{i_hi_{h-1}}}$, $h\ge 2$. Then, a basis of the space $\mJ^2$ is
given by
\begin{align}\label{eqn:basisJ2}
\phi_C(\{x_i\}_{i\in X}), \quad C\in \mR'.
\end{align}
We denote by $\widehat{\B_2}(X,q)$ the quadratic approximation of $\B(X,q)$,
that is the algebra defined by relations $\langle\mJ^2\rangle$. For more details see
\cite[Lemma 2.2]{GG}.

Let $G$ be a finite group. A \textit{principal YD-realization} over
$G$ of $(X,q)$, \cite[Def. 3.2]{AG2}, is a  way to realize this braided vector
space $(\ku X,c^q)$ as a Yetter-Drinfeld module over $G$. Explicitly, it is a
collection
$(\cdot\, , g, (\chi_i)_{i\in X})$ where \begin{itemize}\item $\cdot$ is an
action
of $G$
on $X$,
\item $g:X\to G$ is a function such that $g_{h\cdot i} = hg_
{i}h^{-1}$ and $g_{i}\cdot j=i\rhd j$,
\item the family
$(\chi_i)_{i\in X}$, where $\chi_i:G\to\ku^*$ is a 1-cocycle,
that is
$$\chi_i(ht)=\chi_i(t)\, \chi_{t\cdot i}(h),$$ for all
$i\in X$, $h,t\in G$, satisfying $\chi_i(g_{j})=q_{ji}$.\end{itemize}

If $(\cdot\, , g, (\chi_i)_{i\in X})$ is a principal YD-realization of
$(X,q)$ over $G$ then $\ku X\in  \ydg$ as follows. The action and coaction of
$G$ are determined by:
$$ \delta(x_i)=g_i\ot x_i, \quad h\cdot x_i = \chi_i(h)\, x_{h\cdot i} \quad
i\in X, h\in G.$$

\begin{lema} Assume that for any pair $i, j\in X$, $(i\rhd j)\rhd i=j$, then
\begin{align}\label{chi-q}
    \chi_i(f)\,q_{f\cdot i\rhd  f\cdot j,f\cdot i}  = \chi_j(f)\, q_{i\rhd j,i}
\quad \text{ for any } f\in G,\, i, j\in X.\qed
\end{align}
\end{lema}

\subsection{Nichols algebras over $\s_n$}\label{subsec:nichols}

Let $X=\mO_2^n$ or $X=\mO_4^4$ considered as racks with the map $\rhd$ given by
conjugation. Consider the applications:
\begin{align*}
  \sg&:\s_n\times X\to \ku^*, && (\sigma,i)\mapsto \sg(\sigma), \\
 \chi&: \s_n\times \mO_2^n\to \ku^*, && (\sigma,i)\mapsto
\chi_i(\sigma)=\begin{cases}
  1,  & \mbox{if} \  i=(a,b) \text{ and } \sigma(a)<\sigma(b) \\
  -1, & \mbox{if} \ i=(a,b) \text{ and } \sigma(a)>\sigma(b).
\end{cases}
\end{align*}
We will  deal with the cocycles:
\begin{align*}
 -1&:X\times X\to \ku^*, && (j,i)\mapsto \sg(j)=-1, && i,j\in X;\\
 \chi&: \mO_2^n\times \mO_2^n\to \ku^*, && (j,i)\mapsto \chi_i(j)&& i,j\in\mO_2^n.
\end{align*}

The quadratic approximations of the corresponding Nichols algebras are
\begin{align*}
\widehat{\B_2}(\mO_2^n,-1)=  \,\ku\langle &x_{(lm)}, 1\le l < m \le
n \,\vert\,
x_{(ab)}^2,\, x_{(ab)}x_{(ef)}+x_{(ef)} x_{(ab)}, \\
&x_{(ab)}x_{(bc)}+x_{(bc)} x_{(ac)}+x_{(ac)} x_{(ab)}, \\
& 1\le a < b < c \le n, 1\le e < f \le n, \{a,b\}\cap\{e,f\}=\emptyset
\rangle,
\end{align*}
\begin{align*}
\widehat{\B_2}(\mO_2^n,\chi)=  \,\ku\langle &x_{(lm)}, 1\le l < m
\le n \,\vert\,
x_{(ab)}^2,\, x_{(ab)}x_{(ef)} - x_{(ef)} x_{(ab)},  \\
&x_{(ab)}x_{(bc)} - x_{(bc)} x_{(ac)} - x_{(ac)} x_{(ab)},
\\&x_{(bc)}x_{(ab)} -  x_{(ac)}x_{(bc)} - x_{(ab)}x_{(ac)}, \\
& 1\le a < b < c \le n, 1\le e < f \le n, \{a,b\}\cap\{e,f\}=\emptyset
\rangle,\\
\widehat{\B_2}(\mO_4^4,-1) =\, \ku\langle &x_i, i\in \mO_4^4 \vert
x_i^2,\, x_i
x_{i^{-1}}+x_{i^{-1}} x_i, \\
& x_i x_j+x_k x_i+x_j x_k,  \mbox{if }
ij=ki \ \mbox{and} \ j\neq i\neq k\in
\mO_4^4 \rangle.
\end{align*}

\begin{exa}\label{exa:YD-real} A principal YD-realization of $(\mO_2^n,-1)$ or
$(\mO_2^n,\chi)$, respectively
$(X,q)=(\mO_4^4,-1$), over $\s_n$, respectively $\s_4$,  is given by the
inclusion $X\hookrightarrow \s_n$ and the action $\cdot$ is the
conjugation. The family $\{\chi_i\}$ is determined by the cocycle.
In any case $g$ is injective. For $n=3,4,5$,
this is in fact the only possible realization over $\s_n$.
\end{exa}

\begin{rmk}\label{rem:rrprima}
Notice that all $(\mO_2^n,-1)$, $(\mO_2^n,\chi)$, for any $n$ and
$(\mO_4^4,-1)$ satisfy that $\mR=\mR'$.

When $n=3,4,5$, we have from \cite{AG1, GG}
$\widehat{\B_2}(\mO_2^n,-1)=\B(\mO_2^n,-1)$,
$\widehat{\B_2}(\mO_2^n,\chi)=\B(\mO_2^n,\chi)$ and $\dim
\B(\mO_2^n,-1), \dim \B(\mO_2^n,\chi)<\infty$.
\end{rmk}

\subsection{Pointed Hopf algebras constructed from racks}\label{hopf-ql}

A \emph{quadratic lifting datum} $\mQ=(X,q,G, (\cdot, g,
(\chi_l)_{l\in X}),(\gamma_C)_{C\in \mR'})$, or ql-datum,
\cite[Definition 3.5]{GG}, is a collection consisting of a rack $X$, a
2-cocycle $q$, a finite group $G$, a principal YD-realization
$(\cdot, g, (\chi_l)_{l\in X})$ of $(X,q)$ over $G$ such that
$g_i\neq g_jg_k$, for all $i,j,k\in X$, and a collection
$(\gamma_C)_{C\in \mR'}\in\ku$ satisfying that for each
$C=\{(i_2,i_1),\ldots,(i_n,i_{n-1})\}\in\mR'$, $k\in X$,
\begin{align}
  & \gamma_C=0, && \mbox{if } \ g_{i_2}g_{i_1}=1,\label{eqn:lambda1}\\
  & \gamma_C=q_{k
i_{2}}q_{k i_{1}}\gamma_{k\rhd C}, && \mbox{if } k\rhd C=\{k\rhd
(i_2,i_1),\ldots,k\rhd (i_n,i_{n-1})\}\label{eqn:lambda2}.
\end{align}
To each ql-datum $\mQ$ there is attached a pointed Hopf algebra
$\mH(\mQ)$ generated as an algebra by $\{a_l, H_t : l\in X, \,
t\in G\}$ subject to relations:
\begin{align}
\label{eqn:hq1} H_e &=1, \quad H_tH_s=H_{ts}, && t,s\in G;\\
\label{eqn:hq2} H_t\, a_l &= \chi_l(t)\, a_{t\cdot l}\, H_t, && t\in G, \, l \in
X; \\
\label{eqn:hq3} \phi_C(\{a_l\}_{l\in X}) &= \gamma_C(1-H_{g_ig_j}),
&& C\in \mR', \, (i,j)\in C.
\end{align}
Here $\phi_C$ is as in \eqref{eqn:phiC}
above. The algebra $\mH(\mQ)$ has a structure of pointed Hopf algebra setting
$$\Delta(H_t)=H_t\ot H_t,\;\;  \Delta(a_i)=g_i\ot a_i+a_i\ot 1,\quad t\in
G, i\in X.$$ See \cite{GG} for
further details.

\subsection{Pointed Hopf algebras over $\s_n$}\label{pointe-over-sn}
The following ql-data provide examples of (possibly infinite-dimensional)
pointed Hopf algebras over $\s_n$.
\begin{itemize}
    \item[1.] $\mQ_n^{-1}[t]=(\s_n,\mO_2^n,-1,\cdot,\iota,\{0,\alpha,\beta\})$,
    \item[2.]
    $\mQ_n^{\chi}[\lambda]=(\s_n,\mO_2^n,\chi,\cdot,\iota,\{0,0,\alpha\})$ and
    \item[3.] $\mD[t]=(\s_4,\mO_4^4,-1,\cdot,\iota,\{\alpha,0,\beta\})$,
\end{itemize}
for $\alpha,\beta,\lambda\in\ku$, $t=(\alpha,\beta)$. We will present
explicitly the algebras $\hq$ associated to these data. It follows that
relations \eqref{eqn:hq3} for each $C\in\mR'$ with the same cardinality are
$\s_n$-conjugated. Thus it is enough to consider a single relation
for each $C$ with a given number of elements.
\begin{itemize}
 \item[1.]
$\mH(\mQ_n^{-1}[t])$ is the
algebra presented by generators $\{a_i, H_r :
i\in\mO_2^n,r\in\s_n\}$ and relations:
\begin{align*}
&H_e=1, \quad H_rH_s=H_{rs}, \qquad r,s\in\s_n,\\
&H_j a_i=-a_{jij}H_j, \qquad\qquad\quad i,j\in \mO_2^n, \\
& a_{(12)}^2=0;\\
& a_{(12)} a_{(34)}+a_{(34)} a_{(12)}=\alpha(1-H_{(12)}H_{(34)});\\
& a_{(12)} a_{(23)}+a_{(23)} a_{(13)}+a_{(13)}
a_{(12)}=\beta(1-H_{(12)}H_{(23)}).
\end{align*}
\item[2.]  $\mH(\mQ_n^{\chi}[\lambda])$ is the algebra presented
by generators
$\{a_i, H_r : i\in\mO_2^n,r\in\s_n\}$ and relations:
\begin{align*}
&H_e=1, \quad H_rH_s=H_{rs}, \quad\quad r,s\in\s_n,\\
&H_j a_i=\chi_i(j)\, a_{jij}H_j, \qquad\qquad i,j\in \mO_2^n, \\
& a_{(12)}^2=0;\\
& a_{(12)} a_{(34)}-a_{(34)} a_{(12)}=0;\\
&a_{(12)}a_{(23)}-a_{(23)}a_{(13)}-a_{(13)}a_{(12)}=\alpha(1-H_{(12)}H_{(23)}).
\end{align*}
\item[3.] $\mH(\mD[t])$ is the algebra generated by generators
$\{a_i, H_r :
i\in\mO_4^4,r\in\s_4\}$ and relations:
\begin{align*}
\qquad\qquad& H_e=1, \quad H_rH_s=H_{rs}, \qquad r,s\in\s_n,\\
&H_j a_i=-a_{jij}H_j, \qquad\qquad\quad i\in \mO_4^4,\quad j\in\mO_2^4, \\
& a_{(1234)}^2=\alpha(1-H_{(13)}H_{(24)}),\\
& a_{(1234)}a_{(1432)}+a_{(1432)}a_{(1234)}=0,\\
&a_{(1234)}a_{(1243)}+a_{(1243)}a_{(1423)}+a_{(1423)}a_{(1234)}=\beta(1-H_{
(12)}H_{(13)}).
\end{align*}
\end{itemize}

These Hopf algebras have been defined in  \cite[Def. 3.7]{AG1}, \cite[Def.
3.9]{GG}, \cite[Def. 3.10]{GG} respectively. Each of these algebras $\hq$
satisfy $\gr\hq=\widehat{\B_2}(X,q)\#\ku G$, for $G=\s_n$, $n$ as appropriate
\cite[Propositions 5.4, 5.5, 5.6]{GG}.

\begin{rmk}\label{rem:class res} Classification results:
\begin{itemize}
 \item[(a)]\cite{AHS} $\mH(\mQ_3^{-1}[t])$, $t=(0,0)$ or $t=(0,1)$ are all the
non-trivial finite-dimensional pointed Hopf algebras over $\s_3$.
\item[(b)]\cite{GG} $\mH(\mQ_4^{-1}[t])$, $\mH(\mQ_4^{\chi}[\zeta])$,
$\mH(\mD[t])$,
$t\in\mathbb{P}_\ku^1\cup\{(0,0)\}$, $\zeta\in\{0,1\}$ is a complete list of the
non-trivial finite-dimensional pointed Hopf
algebras over $\s_4$.
\end{itemize}
\end{rmk}

We will classify module categories over the category of representations of any
pointed Hopf algebra over $\s_3$ or $\s_4$, that is, of the algebras listed in Remark \ref{rem:class res}.

\section{Coideal subalgebras of quadratic Nichols algebras}\label{sec:coideal}

A fundamental piece of information to determine simple comodule algebras is the
computation of homogeneous coideal subalgebras inside the Nichols algebra. This
is part of Theorem \ref{mod-over-pointed}. The study of coideal subalgebras is
an active field of research in the theory of Hopf algebras and quantum groups,
see for example \cite{HK}, \cite{HS}, \cite{K} and \cite{KL}.

\medbreak

In this section we present a description of all homogeneous left coideal
subalgebras in the quadratic approximations of the Nichols algebras constructed
from racks.

\medbreak

Fix $n\in\Na$, $X=\{i_1,\dots, i_n\}$ a rack of $n$ elements and $q:X\times X\to\ku^*$ a
2-cocycle. Let $\mR$ be as in \ref{quadratic-app}. Assume that, for any
equivalence class $C$ in $\mR$ and $i,j,k\in X$,
\begin{align}
\label{eqn:ij} &(i,j), (i,k)\in C\Rightarrow j=k \quad \text{and} \quad (i,j),
(k,i)\in C\Rightarrow k=i\rhd j.
\end{align}
Let $G$ be a finite group and let $(\cdot, g, (\chi_i)_{i\in X})$ be a principal
YD-realization of $(X,q)$
over $G$. We shall further assume that
\begin{align}
\label{eqn:gmono} &g \quad \text{ is injective and  } \mR=\mR'.
\end{align}

For each subset $Y\subseteq X$, $Y=\{i_{j_1},\dots,i_{j_r}\}\subseteq X$, denote
by
$\mK_Y$ the subalgebra of
$\widehat{\B_2}(X,q)\#\ku 1$ generated by $ x_{j_1},\dots, x_{j_r}$.
Set $\mH=\widehat{\B_2}(X,q)\#\ku G$.

\begin{prop}\label{pro:stabilizers}
For each set $Y=\{i_{j_1},\dots,i_{j_r}\}\subseteq X$ the algebra $\mK_Y$ is an
homogeneous coideal subalgebra of $\mH$. For each such selection, if
$S=\{g_{i_{j_1}},\dots,g_{i_{j_r}}\}$ then
$$\stab \mK_Y=S_Y^G=\{h\in G: hS_Yh^{-1}= S_Y\}.$$

Moreover, if $\mK$ is a homogeneous coideal subalgebra of $\mH$ generated in degree
one, then there exists a unique $Y\subseteq X$ such that
$$
\mK=\mK_Y.$$
In particular, the set of homogeneous coideal subalgebras of $\mH$ generated in degree
one inside $\widehat{\B_2}(X,q)\#\ku 1$
is in bijective correspondence with the set $2^X$ of parts of $X$.

\end{prop}
\begin{proof}
It is clear that $\mK=\mK_Y$ is a homogeneous
coideal subalgebra. Now, to describe $\stab \mK$ it is enough to
compute the stabilizer of the vector space $\ku\{ x_{j_1},\dots,
x_{j_r}\}$. But $h\cdot x_{j_k}=\chi_{j_k}(h)x_{h\cdot j_k}$,
$k=1,\dots,r$ and $x_{h\cdot j_k}\in \{ x_{j_1},\dots, x_{j_r}\}$ if
and only if $h\cdot j_k\in\{j_1,\dots,j_r\}$, if and only if
$g_{h\cdot j_k}=g_{j_l}$ for some $l=1,\dots,r$. And the first part of the proposition
follows since $g_{h\cdot j_k}=hg_{j_k}h^{-1}$.

Now, let $\mK$ be a homogeneous coideal subalgebra of $\mH$ generated in degree
one. If $\mK=\ku$ the result is trivial, so let us assume that $\mK\neq \ku$. Since
$\mK$ is homogeneous then $\mK(1)\neq 0$.
Let $0\neq y=\sum_i\lambda_i x_i\in\mK(1)$, then
$$
\Delta(y)=y\ot 1 +\sum_i \lambda_i H_{g_i} \ot
x_i\Rightarrow \sum_i \lambda_i H_{g_i} \ot
x_i\in \mH_0 \ot \mK(1).
$$
Let $\sum_i \lambda_i H_{g_i} \ot
x_i=\sum_{t\in G}H_t\ot \kappa_t$, $\kappa_t=\sum_{j\in
X}\eta_{tj}x_j\in\mK(1)$, $\eta_{tj}\in\ku$, $\forall\,t,j$.

As $H_t=H_{g_j}$ if and only if $t=g_j$ and $g_i=g_j$ if and only if $i=j$,
for every $i,j\in X,\,t\in G$, \eqref{eqn:gmono}, then $\eta_{tk}=0$ if
$t\neq g_k$ for some $k\in X$. Let us denote $\eta_{ij}=\eta_{g_ij}$, thus,
\begin{align*}
&\sum_i \lambda_i H_{g_i} \ot
x_i=\sum_{i,j}\eta_{ij}H_{g_i}\ot x_j.
\end{align*}
Therefore, $\lambda_i\neq 0$ implies $\eta_{ij}=\delta_{i,j}\lambda_i$ and so
$\kappa_i=x_i$. Thus, $\{x_i\,|\,\lambda_i\neq
0\}\subset \mK$, $\mK(1)=\bigoplus\limits_{x_i\in\mK(1)}\ku\{x_i\}$ and therefore if $Y=\{i\in X\,:\,x_i\in\mK(1)\}$ then $\mK=\mK_Y$.
Finally, if
$Y\neq Y'$ then it follows from the injectivity of $g$ that $\mK_Y\ncong
\mK_{Y'}$ as coideal subalgebras.
\end{proof}

The next general lemma will be useful in  \ref{subsec:subcoideales de s3} to prove that certain subalgebras are generated in degree one. Given a rack $X$,
let us recall the notion of derivations $\delta_i$ associated to every element of the canonical basis
$\{e_i\}_{i\in X}$. If $\{e^i\}_{i\in X}$ denotes the dual basis of this basis, then $\delta_i=(\id\ot e^i)\Delta$. If $i\in X$ we denote by $X_i$ the
set $X\setminus\{i\}$ and thus $\ku X_i=\ku\{x_j\,|\,j\in X_i\}$. Let us assume, furthermore, that
\begin{equation}\label{eqn:qii=-1}
q_{ii}=-1, \quad \text{ for all } \,i\in X.
\end{equation}
This condition is satisfied, for example, if $\dim\widehat{\B_2}(X,q)<\infty$ or
$X$ is such that $i\rhd i=i$ for all  $\,i$, by \eqref{eqn:gmono}.
\begin{lema}\label{lem:derivaciones}
Let $\mK\subset \widehat{\B_2}(X,q)\#\ku 1$ be a homogeneous coideal subalgebra of $\mH$. Let $i\in X$ and let us assume that there is $\omega\in\mK$ such that $\delta_i(\omega)\neq0$. Then $x_i\in\mK(1)$.
\end{lema}

\pf
Let $\mK=\bigoplus_s \mK(s)$, $\omega\in T(\ku X)$ and $i\in X$. In $\mH$,
$$
\omega=\alpha_i(\omega)+\beta_i(\omega)x_i, \qquad \alpha_i(\omega),
\beta_i(\omega)\in \mK_{X_i}.
$$
It suffices to see this for a homogeneous monomial $\omega$. We see it by induction in $\ell=\ell(\omega)\in\N$ such that $\omega\in T^\ell(\ku X)$. If $\ell=0$, or
$\ell=1$ this is clear. Let us assume it holds for $\ell=n-1$, for some $n\in\N$. If  $\ell(\omega)=n$ and $\omega=x_{j_1}\dots x_{j_n}$, two possibilities hold, that is $j_1\neq i$ or $j_1= i$. In the first case, let
$\omega'=x_{j_2}\dots x_{j_n}$. Thus, $\ell(\omega')\leq n-1$ and therefore there exist $\alpha_i(\omega'), \beta_i(\omega')\in \mK_{X_i}$ such that
$\omega'=\alpha_i(\omega')+\beta_i(\omega')x_i$. As
$x_{j_1}\alpha_i(\omega')$, $x_{j_1}\beta_i(\omega')\in \mK_{X_i}$ the claim follows in this case.

In the second case, let $j=j_2$ and let us note that $j\neq i$, by \eqref{eqn:qii=-1}.
By \eqref{eqn:gmono}, we can consider the relation
$$
x_ix_j=q_{ij}x_{i\rhd j}x_i-q_{ij}q_{i\rhd j\,i}x_jx_{i\rhd j}.
$$
Thus, if $\omega''=x_{j_3}\dots x_{j_n}$, $\omega=q_{ij}x_{i\rhd
j}x_i\omega''-q_{ij}q_{i\rhd j\,i}x_jx_{i\rhd j}\omega''$ and both members of this sum belong to $\mK_{X_i}+\mK_{X_i}x_i$ because of the previous case and thus the claim follows.

\medbreak

Let $\pi:\bigoplus\limits_{s=0}^m \mH(s)\ot \mK(m-s)\to \mH(m-1)\ot \mK(1)$ be the canonical linear projection. Let $\omega\in T(\ku X)$, $i\in X$ and
$\alpha_i(\omega)$, $\beta_i(\omega)$ as above. Then,
$$\pi\Delta(\omega)\in
\beta_i(\omega)\ot x_i+\bigoplus\limits_{j\neq i}\mH\ot x_j.$$
Notice that $\delta_i(\omega)=\beta_i(\omega)$, and therefore if
$\delta_i(\omega)\neq 0$ it follows that  $x_i\in\mK(1)$ using
\eqref{eqn:gmono} as in the proof of Proposition \ref{pro:stabilizers}.
\epf

In this part we shall assume that $X$ is one of the racks $\mO_2^n$, $n\in\N$,
or $\mO_4^4$, $q$ one of the cocycles in \ref{subsec:nichols}. Notice that
\eqref{eqn:ij} is satisfied in these cases. Using the previous results we shall
describe explicitly all connected
homogeneous coideal subalgebras of the bosonization of the quadratic
approximations to Nichols algebras described in \ref{subsec:nichols}.

\medbreak

We first introduce the following notation. Let $Y\subset X$ be a
subset and define
\begin{align*}
 \mR^Y_1&=\{C\in \mR: C\subseteq Y\times Y\},\\
 \mR^Y_2&=\{C\in \mR:\, \mid C\cap Y\times Y\mid= 1\}, \ \text{and} \\
 \mR^Y_3&=\{C\in \mR: C\cap Y\times Y=\emptyset\}.
\end{align*}

\begin{rmk}
For the ql-data in \ref{pointe-over-sn}, $\mR=\mR^Y_1 \cup \mR^Y_2\cup \mR^Y_3$, for any subset $Y$.
Moreover, if $f\in \stab\mK_Y$ then $f \cdot \mR^Y_s \subseteq \mR^Y_s$ for any
$s=1,2,3$. Also, \eqref{eqn:qii=-1} holds.

\end{rmk}

\begin{defi}\label{fed:LY}
Take the free associative
algebra $\mT$ in the variables $\{T_l\}_{l\in Y}$. According to this, we set
$\vartheta_{C,Y}(\{T_l\}_{l\in Y})$ in $\mT$ as
\begin{equation}\label{eqn:psiC}
  \vartheta_{C,Y}(\{T_l\}_{l\in Y})=\begin{cases}
    \phi_C(\{T_l\}_{l\in X}), &\text{ if } C\in \mR^Y_1; \\
    T_iT_jT_i-q_{i\rhd j,i}\; T_jT_iT_j, & \text{ if } C\in \mR^Y_2, (i,j)\in C
\cap Y\times Y; \\
    0, & \text{ if } C\in \mR^Y_3.
  \end{cases}
\end{equation}
We define the algebra $\mL_Y$ as follows
\begin{equation}\label{eqn:relaciones}
\mL_Y= \ku\langle \{y_i\}_{i\in Y} \rangle /\langle
\vartheta_{C,Y}(\{y_l\}_{l\in Y}): C\in\mR\rangle.
\end{equation}
\end{defi}

Notice that, if $Y=X$ then $\mL_X\cong\B(X,q)$. For simplicity we shall
sometimes denote $\vartheta_C=\vartheta_{C,Y}$.

\medbreak

We now take $\B$ one of the quadratic (Nichols) algebras
$\widehat{\B_2}(\mO_2^n,-1)$,
$\widehat{\B_2}(\mO_2^n,\chi)$, or $\B(\mO_4^4,-1)$. Accordingly, let
$X=\mO_2^n$, $q=-1,\chi$ or $(X,q)=(\mO_4^4,-1)$. Consider a YD-realization for
$(X,q)$ such that \eqref{eqn:gmono} is satisfied (for instance, the ones in
Example \ref{exa:YD-real}). Set $\mH=\B\#\ku G$.

\begin{teo}\label{teo:relations}
Let $Y\subset X$. $\mL_Y$ is an $\mH$-comodule algebra with coaction
$$
\delta(y_i)=g_i\ot y_i+x_i\ot 1, \quad i\in Y.
$$
The application $y_i\mapsto x_i$, $i\in Y$ defines an epimorphism of
$\mH$-comodule algebras $\mL_Y\twoheadrightarrow\mK_Y$. Moreover, if $n=3$, it is
an isomorphism and $\mL_Y\cong\mK_Y$.
\end{teo}
\pf
The relations that define $\mL_Y$ are satisfied in $\B$. In fact, it suffices to
check this in the case $C\in \mR^Y_2$ since in the other ones $\vartheta_C=0$ or
$\vartheta_C=\phi_C$, and $\phi_C=0$ in $\B$, see \eqref{eqn:basisJ2}. Now, if
$C\in \mR^Y_2$ and $(i,j)\in C \cap Y\times Y$, let $k=i\rhd j$. By the
definition of $\mR^Y_2$, we have that necessarily $k\neq i,j$. Then, if we
multiply the relation $x_ix_j-q_{ij}x_{i\rhd j}x_i+q_{ij}q_{i\rhd
j\,i}x_jx_{i\rhd j}=0$ by $x_i$ on the right and apply this relations to the
outcome, we get
    \begin{align*}
0&=x_ix_jx_i+q_{ij}q_{i\rhd j\,i}\, x_jx_{i\rhd j}x_i=x_ix_jx_i+q_{i\rhd
j\,i}x_j(x_ix_j+q_{ij}q_{i\rhd j\,i}\, x_jx_{i\rhd j})\\
&=x_ix_jx_i+q_{i\rhd j\,i}\, x_jx_ix_j.
    \end{align*}
Thus, we have an algebra projection $\pi:\mL_Y\twoheadrightarrow\mK_Y$. It is
straightforward to see that, for every $C\in\mR$,
$$
\delta(\vartheta_{C,Y}(\{y_l\}_{l\in Y}))=\vartheta_{C,Y}(\{x_l\}_{l\in Y})\ot
1+g_{C,Y}\ot \vartheta_{C,Y}(\{y_l\}_{l\in Y}),
$$
where
\begin{align*}
  g_{C,Y}=\begin{cases}
    g_ig_j &\text{if} \quad C\in \mR^Y_1, (i,j)\in C, \\
    g_ig_jg_i & \text{if} \quad C\in \mR^Y_2, (i,j)\in C \cap Y\times Y, \\
    0, & \text{if} \quad C\in \mR^Y_3.
  \end{cases}
\end{align*}
Therefore, $\delta$ provides $\mL_Y$ with a structure of $\mH$-comodule in such a
way that $\pi$ becomes a homomorphism.

\medbreak

We analyze now the particular case $n=3$. If $|Y|=1$, the result is clear. Let
us suppose then that $Y=\{i,j\}\subset\mO_2^3$. Notice that the map $\pi$ is
homogeneous. Then, if  $\gamma\in\ker(\pi)$, $\pi(\gamma)=0$ in
$\B(\mO_2^3,-1)$. By \eqref{eqn:basisJ2}, we have that necessarily
$\deg\gamma\geq 3$. Now, if $\deg\gamma=3$,
$$
\gamma=\alpha y_iy_jy_i+\beta y_jy_iy_j=(\alpha+\beta)y_jy_iy_j.
$$
for $\alpha,\beta\in\ku$. Then, $\pi(\gamma)=0$ implies that $\alpha=-\beta$ and
$\gamma=0$. Finally, we can see that there are no elements $\gamma\in\mL_Y$ with
$\deg\gamma\geq4$. In fact, an element $\gamma$ with $\deg\gamma=4$ would be of
the form
$$
\gamma=\alpha y_iy_jy_iy_j+\beta y_jy_iy_jy_i=\alpha y_iy_iy_jy_i+\beta
y_jy_jy_iy_j=0.
$$
This also shows that there are no elements of greater degree. Therefore,
$\mL_Y=\mK_Y$.
\epf

\begin{rmk}
If $n\neq 3$, then in general $\mL_Y\neq\mK_Y$. In fact, if $n=4$, $q=-1$ and we
take $Y=\{(13), (23), (34)\}\subseteq \mO_2^4$, then we have
\begin{align*}
\mL_Y\cong \ku\langle x,y,z\,:\,  x^2, y^2, z^2,  xyx-yxy, yzy-zyz,
xzx-zxz\rangle.
\end{align*}
Now, in the subalgebra of $\B(\mO_2^4,-1)$ generated by $x=x_{(23)}$,
$y=x_{(34)}$, $z=x_{(13)}$ we have the relation
\begin{align*}
(xyz)^2&=x_{(23)}x_{(34)}x_{(13)}x_{(23)}x_{(34)}x_{(13)}\\
&=-x_{(23)}x_{(34)}(x_{(23)}x_{(12)}+x_{(12)}x_{(13)})x_{(34)}x_{(13)}\\
&=x_{(23)}x_{(34)}x_{(23)}x_{(34)}x_{(12)}x_{(13)}\\
&\quad + x_{(23)}x_{(12)}x_{(34)}x_{(13)}x_{(34)}x_{(13)}\\
&=x_{(23)}x_{(23)}x_{(34)}x_{(23)}x_{(12)}x_{(13)}\\
&\quad + x_{(23)}x_{(12)}x_{(34)}x_{(34)}x_{(13)}x_{(34)}\\
&=0.
\end{align*}
But $(xyz)^2\neq 0$ in $\mL_Y$. We prove this by using the computer program
\cite{GAP} together with the package  \cite{GBNP}.
See Proposition \ref{pro:coideal transp} (6), for a description of $\mK_Y$ in
this case.
\end{rmk}

\subsection{Coideal subalgebras of Hopf algebras over $\s_n$}\label{subsec:subcoideales de s3}

Set $n=3,4$, $\B$ a finite-dimensional Nichols algebra over $\s_n$, $\mH=\B\#\ku\s_n$. Recall that these Nichols algebras coincide with their
quadratic approximations. We will describe all
the coideal subalgebras of $\mH$. We will also calculate their
stabilizer subgroups.

We start out by proving that in this case these coideal subalgebras are generated in degree one.

\begin{teo}\label{teo:degree1}
 Let $\mK$ be a homogeneous left coideal
subalgebra of $\mH$. Then $\mK$ is generated in degree one. In particular,
$\mK=\mK_Y$ for a unique $Y\subseteq X$.
\end{teo}
\pf
We will se that, given $\omega\in\mK$, $\omega\in\langle x_i\,:\,\delta_i\omega\neq 0\rangle$. Then, by Lemma \ref{lem:derivaciones}, it will follow that $\omega\in\langle\mK(1)\rangle$. Let $I=\{i\in X\,:\,\delta_i\omega=0\}$ and let us assume $I$ has $m$ elements. We will see this case by case, for $m=0,\dots, 6$.

\medbreak

Cases $m=0$ (that is, for all $i\in X$ $x_i\in \mK(1)$) and $m=6$, or $m=n$ in general, (since in this case $\omega=0$, see \cite[Section 6]{AG1}) are clear.
Case $m=1$ is Lemma \ref{lem:derivaciones}, which also holds for any $n\in\N$.

\medbreak

Let us see case $m=2$, for any $n\in \Na$. Let $I=\{i,p\}$. We know that there is an expression of $\omega$ without, v.g., $x_i$. Let us see that we can write $\omega$ without $x_i$ nor $x_p$. Let $j\in X$ such that $p\rhd j=i$. Using the relations as in Lemma \ref{lem:derivaciones}, and using  that $x_lx_rx_l=-q_{l\rhd r l} x_rx_lx_r$ and that $x_rx_lx_rx_l=0$ for all $l,r\in X$, we can assume that $\omega$ can be written as
$$
\omega=\gamma^0+\gamma^1x_p+\gamma^{2}x_px_j
$$
with $\gamma^0$, $\gamma^1$, $\gamma^{2}$, such that they do not contain factors $x_{i}$, $x_p$ in their expressions. Let us see this in detail. We can assume that $\omega\in T^\ell(\ku X)$ is a homogeneous monomial. For each appearance of a factor $x_px_l$, with $l\neq j$ we change it by $q_{p l}\; x_lx_{p\rhd l}+ q_{pl}q_{p\rhd l p}\; x_{p\rhd l}x_p$. That is, we change for an expression in which $x_p$ is located more to the right and an expression that does not contain $x_i$ nor $x_p$ (in the place where we had an $x_p$). If we have a factor of the form $x_px_j$ we move it to the right, until we get to $x_px_jx_p$, but we can change this expression by $-q_{p\rhd j p}\; x_jx_px_j$.

Now, $0=\delta_p\omega=\gamma^1g_p+\gamma^{2}g_px_j=(\gamma^1+q_{pj}\gamma^{2}x_i)g_p$, and therefore we have
$$
\omega=\gamma^0+\gamma^1x_p+q_{pj}\gamma^{2}x_ix_p+q_{pj}q_{ip}\gamma^{2}x_jx_i=\gamma^0+q_{pj}q_{ip}\gamma^{2}x_jx_i.
$$
But then $0=\delta_i\omega=q_{pj}q_{ip}\gamma^{2}x_jg_i$ and therefore $\omega=\gamma^0$ can be written without $x_i, x_p$.

\medbreak

This finishes the case $n=3$, since in this case $|X|=3$. We now fix $n=4$, to deal with the cases $m=3,4,5$.

\medbreak

Let us see case $m=3$. Fix $I=\{i_1,i_2,p\}$. There are three possibilities
\begin{align}
\label{caso1} I&=\{i,j,i\rhd j\};\\
\label{caso2} I&=\{i,j,k\}, && \text{such that } i\rhd k=k\text{ or }j\rhd k=k;\\
\label{caso3} I&=\{i,j,l\}, && \text{the remaining case}.
\end{align}
Set $j_1,j_2\in X$ be such that $p\rhd j_s=i_s$, $s=1,2$. We can assume that $\omega$ is written without $x_{i_s}$, $s=1,2$. Notice that not always $j_1,j_2$ exist. For instance, in \eqref{caso1} there are no $j_1$ nor $j_2$ and in \eqref{caso2} $j_1$ or $j_2$ does not exist. We analyze the three cases separately.

In \eqref{caso1}, as there are no $j_1$, $j_2$, we can write $\omega$ in the form $\omega=\gamma^0+\gamma^1x_p
$ with $\gamma^0,\gamma^1$ without factors $x_j$, $j\in I$. But from $\delta_p\omega=0$ it follows that $\omega=\gamma^0$ and therefore we can write $\omega$ without factors $x_j$, $j\in I$.

Case \eqref{caso2} is similar. Assume, for example, that $i_2\rhd p=p$. Thus, we have no $j_2$. Accordingly, we can assume that $\omega$ is of the form
\begin{align*}
\omega&=\gamma^0+\gamma^1x_p+\gamma^2x_px_{j_1}=\gamma^0+\gamma^1x_p+q_{p j_1}\gamma^2x_{i_1}x_p+q_{p j_1}q_{j_1 i_1}\gamma^2x_{j_1}x_{i_1}
\end{align*}
with $\gamma^1, \gamma^2, \gamma^3$ without factors $x_j$, $j\in I$. Now,
$
0=\delta_p\omega=(\gamma^1p+q_{p j_1}\gamma^2x_{i_1})g_p
$
and thus
$
\omega=\gamma^0+q_{p j_1}q_{j_1 i_1}\gamma^2x_{j_1}x_{i_1}
$
but as $\delta_{i_1}\omega=0$, it follows $\omega=\gamma^0$ and therefore $\omega$ is written without factors $x_j$, $j\in I$.

It remains to see \eqref{caso3}. The existence of $j_1,j_2$ makes this case more subtle than the previous ones. Let us analyze the set $I=\{i_1,i_2,p\}$. We have that $k=i_1\rhd i_2=i_2\rhd i_1\notin I$, but, moreover, we have that $X=\{i_1,i_2,p,k,j_2,j_1\}$. In fact, we cannot have $i_1\rhd i_2=j_1$ since this implies $i_2=p$ and neither $i_1\rhd i_2=j_2$ because this implies $i_1=p$. Moreover, we have that $i_2\rhd j_1=j_1$, and therefore $x_{i_2}x_{j_1}=q_{i_2j_1}x_{j_1}x_{i_2}$. Set
\begin{align*}
& a=x_{p}, && b=x_{j_1}, && c=x_{j_2},\\
& d=x_{i_1}, && e=x_{i_2}, && f=x_{k}.
\end{align*}
Now we analyze which are the longest words we can write with the ``conflictive'' factors $a$, $b$ and $c$, starting with $a$. Recall that $aba=\pm bab$, and $abb=0$. Starting with $ab$, we can preliminary form the words $abca$ and $abcb$. Now, $abcac=\pm babca$, and thus we discard it. Consider $abcb$. $abcabc=0$, so we are left with $abcaba$. As $abcabab=0$, we reach to $abcabac$. As $abcabaca=abcabacb=0$, we keep this word. In the case of $abcb$, arguing similarly, we reach to $abcbacb$. If we start with $acb$, as $acbc=\pm abcb$, we consider those words starting with $acba$. The longest one is $acbacab$, but this is $\pm abcbacb$. So the longest word we can form not considered before is $acbac$.

In consequence, we can assume there exist $\gamma^i\in\mK$, $i=0,\dots, 15$ without factors $x_j$, $j\in I$, such that $\omega$ is of the form
\begin{align*}
\omega=&\gamma^0+\gamma^1a+\gamma^2ab+\gamma^3abc+\gamma^4abca+\gamma^5abcab+\gamma^6abcaba+\gamma^7abcabac\\
&+\gamma^8abcb+\gamma^9abcba+\gamma^{10}abcbac+\gamma^{11}abcbacb+\gamma^{12}ac+\gamma^{13}acb\\
&+\gamma^{14}acba+\gamma^{15}acbac.
\end{align*}
Using the relations and the fact that $\delta_s\omega=0$, $s=p,i_1,i_2$ we will show that we can write $\omega$ without factors $x_s$, $s=p,i_1,i_2$. When using the relations, by abuse of notation, we will omit the scalars $q_{..}$ that may appear, including them in the (new) factors $\gamma^i$. We will denote by $\gamma^{i'}$, $\gamma^{i''}$, $\gamma^{i'''}\in\mK$ to some of these scalar multiple of the factors $\gamma^i$, $i=0,\dots, 15$, when needed.

As $\delta_p\omega=0$, we can re-write $\omega$ as
\begin{align*}
\omega=&\gamma^0+\gamma^2bd+\gamma^3bdc+\gamma^{3'}dce+\gamma^5bdcbd+\gamma^{5'}dcebd+\gamma^7abcabce\\
&+\gamma^8bdcb+\gamma^{8'}dceb+\gamma^{8''}debd+\gamma^{10}abcbce+\gamma^{11}abcbebd+\gamma^{11'}abcbceb\\
&+\gamma^{12}ce+\gamma^{13}ebd+\gamma^{13'}ceb+\gamma^{15}acbea+\gamma^{15'}acbce.
\end{align*}
Now, using that $\delta_{i_1}\omega=0$, and the relations $dc=\pm cd$, $be=\pm eb$, $bcb=\pm cbc$ and $abcabc=bcbc=0$, we see that
\begin{align*}
\omega=&\gamma^0+\gamma^2bd+\gamma^3bcd+\gamma^{3'}cde+\gamma^5bdcbd+\gamma^{5'}dcebd\\
&+\gamma^7abcdeae+\gamma^8bdcb+\gamma^{8'}dcbe+\gamma^{8''}dbed+\gamma^{11}abcbeda\\
&+\gamma^{12}ce+\gamma^{13}bed+\gamma^{13'}cbe+\gamma^{15}acbce+\gamma^{15'}edaea.
\end{align*}
Using that  $\delta_{i_{2}}\omega=0$ together with the relations, we get to
\begin{align*}
\omega=&\gamma^0+\gamma^2bd+\gamma^3bcd+\gamma^5bcbad+\gamma^{5'}cbafe+\gamma^{5''}cbaed\\
&+\gamma^8bcad+\gamma^{8'}bcba+\gamma^{8''}baed+\gamma^{8'''}bafe+\gamma^{11}abcbebd\\
&+\gamma^{11'}abcbeab+\gamma^{11''}abcbfea+\gamma^{11'''}abcbfac+\gamma^{13}bed+\gamma^{13'}bfe.
\end{align*}
Using now that  $\delta_{i_{1}}\omega=0$,
\begin{align*}
\omega=&\gamma^0+\gamma^5cbafe+\gamma^8bcba+\gamma^{8'}bafe\\
&+\gamma^{11}abcbacb+\gamma^{11'}abcbfce+\gamma^{11''}abcbfac+\gamma^{13}bfe.
\end{align*}
Using again that $\delta_{i_{2}}\omega=0$,
\begin{align*}
\omega&=\gamma^0+\gamma^8bcba+\gamma^{11}abcbacb+\gamma^{11'}abcbafc\\
&=\beta^0+\beta^1a+\beta^2abcbacb+\beta^3abcbabf
\end{align*}
for $\beta^i\in\mK$, $i=0,\dots, 3$ without factors $x_j$, $j\in I$. Using that $\delta_p\omega=0$,
\begin{align*}
\omega=\beta^0+\beta^2dedaeda+\beta^3dedadaf=\beta^0,
\end{align*}
since $edaeda=dada=0$. That is, we can write $\omega$ without factors $x_j$, $j\in I$.

\medbreak

In the case $m=4$, we look at the different subsets of three elements of $I$. If we have a subset of three elements that corresponds to the case \eqref{caso3} it follows that $\omega$ can be written without the factors $x_j$ with $j$ in that subset, and then $\omega$ is in an algebra isomorphic to $\B(\mO_2^3,-1)$, for which we have already proved the result. If we have a subset as in the case \eqref{caso1} when we add to this subset a fourth element we obtain another subset as in the case \eqref{caso3}. If our subset corresponds to the case \eqref{caso2}, in order to get to a case different from \eqref{caso3}, we necessarily have to add a fourth element such that $I$ is
\begin{align*}
&I=\{i,j,k,l\}, && \text{with } i\rhd k=k\text{ and }j\rhd l=l.
\end{align*}
We analyze this case. If $p\in I$ is fixed and $\omega$ is written without factors $x_j$, $j\in I\setminus\{p\}=\{i_1,i_2,i_3\}$, notice that if $p\rhd i_3=i_3$ there is no other $j_3$ such that $p\rhd j_3$ and, moreover, if $j_1,j_2$ are such that $p\rhd j_s=i_s$, $s=1,2$, then $x_{j_1}x_{j_2}=\pm x_{j_2}x_{j_1}$. Therefore, we can assume that there are $\gamma^i$, $i=0,\dots, 4$ such that they do not contain factors $x_j$, $j\in I$ and such that $\omega$ can be written as
\begin{align*}
\omega&=\gamma^0+\gamma^1x_p+\gamma^2x_px_{j_1}+\gamma^3x_px_{j_1}x_{j_2}+\gamma^4x_px_{j_1}x_{j_2}x_{p}\\
&\overset{\delta_p\omega=0}{=} \gamma^0+\gamma^2x_{j_1}x_{i_1}+\gamma^3x_{j_1}x_{i_1}x_{j_2}+\gamma^{3'}x_{i_1}x_{j_2}x_{i_2}+\gamma^4x_{i_1}x_{i_2}x_p\\
&\overset{\delta_p\omega=0}{=}\gamma^0+\gamma^2x_{j_1}x_{i_1}+\gamma^3x_{j_1}x_{i_1}x_{j_2}+\gamma^{3'}x_{i_1}x_{j_2}x_{i_2}\\
&\overset{\delta_{i_2}\omega=0}{=}\gamma^0+\gamma^2x_{j_1}x_{i_1}+\gamma^3x_{j_1}x_{i_1}x_{j_2}\\
&\overset{\delta_{i_1}\omega=0}{=}\gamma^0+\gamma^3x_{j_1}x_{j_2}x_{i_2}\\
&\overset{\delta_{i_2}\omega=0}{=}\gamma^0.
\end{align*}
Then, we can write $\omega$ without $x_j$, for $j\in I$. In the case $m=5$, $\omega$ necessarily belongs to an algebra isomorphic to $\B(\mO_2^3,-1)$.
\epf

Now we apply Theorems \ref{teo:relations} and \ref{teo:degree1} to calculate the coideal subalgebras and stabilizer subgroups of $\mH=\B(\mO_2^3)\#\ku\s_3$.

\begin{cor}\label{cor:coideal s3}
The following are all the  proper homogeneous left coideal
subalgebras of $\B(\mO_2^3,-1)\#\ku\s_3$:
\begin{enumerate}
  \item $\mK_i=\langle x_i\rangle\cong \ku[x]/\langle x^2\rangle$,
$i\in\mO_2^3$;
  \item $\mK_{i,j}=\langle x_i,x_j\rangle\cong \ku\langle
x,y\rangle /\langle
  x^2,y^2,xyx-yxy\rangle$, $i,j\in\mO_2^3$.
\end{enumerate}
The non trivial stabilizer subgroups of $\s_3$ are, on each
case
\begin{enumerate}
  \item $\stab \mK_i=\Z_2\cong \langle i\rangle\subset \s_3$;
  \item $\stab \mK_{i,j}=\Z_2\cong \langle k\rangle\subset \s_3$, $k\neq
  i,j$.\qed
\end{enumerate}
\end{cor}

Next, we use the computer program \cite{GAP}, together with the package
\cite{GBNP}, to compute the coideal subalgebras of the finite-dimensional
Nichols algebras over $\s_4$ associated to the rack of transpositions $\mO_2^4$. In the same way, the coideal
subalgebras of the Nichols algebra $\B(\mO_4^4,-1)$ associated to the rack of
4-cycles can be computed. The presentation of these algebras may not be minimal, in the sense that there
may be redundant relations. Moreover, in the general case, non-redundant
relations in a coideal subalgebra $\mK$ may become redundant when computing the
bosonization with a subgroup $F\leq\stab \mK$.

\medbreak

First, we need to establish some notation and
conventions. Let $\ku\langle x,y,z\rangle$ be the free algebra in
the variables $x,y,z$. We set the ideals
\begin{align*}
R^\pm(x,y,z)&=\langle x^2,y^2,z^2,xy+yz\pm zx\rangle\subset
\ku\langle x,y,z\rangle.
\end{align*}
Set $\B^+_4=\B(\mO_2^4,-1)$, $\B^-_4=\B(\mO_2^4,\chi).$
Recall that $Y$ stands for a subset of $\mO_2^4$.

\begin{prop}\label{pro:coideal transp}
Let $\varepsilon=\pm$.Any homogeneous proper coideal subalgebra $\mK^\veps$ 
of $\B_4^\veps\#\ku 1$ is isomorphic to one of the algebras in
the following list:

\noindent $\dim\mK^\veps(1)=1$,
\begin{enumerate}
  \item[(1)] $Y=\{i\}$, $\mK^\veps= \ku[x]/\langle x^2\rangle,$ and $\dim\mK^\veps=2$.
\end{enumerate}

\noindent$\dim\mK^\veps(1)=2$,
\begin{enumerate}
  \item[(2)] $Y=\{i,j\}$, $i\rhd j=j$,

$\mK^\veps= \ku\langle
x,z\rangle /\langle  x^2,z^2,xz+\veps zx\rangle,$ and $\dim\mK^\veps=4$.
\item[(3)] $Y=\{i,j\}$, $i\rhd j\neq j$,

$\mK^\veps= \ku\langle
x,y\rangle /\langle
  x^2,y^2,xyx-\veps yxy\rangle,$ and  $\dim\mK^\veps=6$.
\end{enumerate}

\noindent $\dim\mK^\veps(1)=3$:
\begin{enumerate}
\item[(4)] $Y=\{i,j,k\}$, $i\rhd j=k$,

$\mK^\veps= \ku\langle
x,y,z\rangle /\langle
  R^\veps(x,y,z)\rangle,$
 and  $\dim\mK^\veps=12$.
\item[(5)] $Y=\{i,j,k\}$, $i\rhd j \neq j,k$, $i\rhd k=k$

$\mK^\veps_{i,j,k}:= \ku\langle
x,y,z\rangle /\langle
  x^2,y^2,z^2,xyx-\veps yxy,zyz-\veps yzy,xz+\veps zx\rangle,$ and
$\dim\mK^\veps=24$.
\item[(6)] $Y=\{i,j,k\}$, $i\rhd j, j\rhd k, i\rhd k\notin\{i,j,k\}$,
\begin{align*}
\mK_Y^\veps=\ku\langle x,y,z\,:\, &  x^2,  y^2,  z^2, \\
 & yxy - \veps xyx, zxz -\veps xzx,  zyz - \veps yzy, \\
& zxyz + yzxy + xyzx,  zyxz + yxzy + xzyx \\
 & zxyxzx + \veps yzxyxz, zxyxzy + \veps xzxyxz \rangle,
\end{align*}
and  $\dim\mK^\veps=48$.
\end{enumerate}

\noindent $\dim\mK^\veps(1)=4$:
\begin{enumerate}
\item[(7)] $Y=\{i,j,k,l\}$, $i\rhd j= k$, $i\rhd l=l$;
\begin{align*}
\mK_Y^\veps=\ku\langle x,y,z,w\,:\, &  x^2 , y^2, z^2, w^2, \\
& zx +\veps yz +\veps xy ,  zy + yx +\veps xz ,  wz +\veps zw,  \\
&yxy - \veps xyx ,  wxw -\veps xwx , wyw -\veps ywy, \\
& wyx +\veps wxz -\veps zwy ,  wyz + wxy - zwx \\
& wxyz - zwxz, wxzw + xwxz , \\
&  wxyw + ywxy + xywx, wxyxz -\veps zwxyx ,\\
& wxyxwx +\veps ywxyxw, wxyxwy +\veps xwxyxw \rangle,
\end{align*}
and  $\dim\mK^\veps=96$.
\item[(8)] $Y=\{i,j,k,l\}$, $i\rhd j\neq j, k$, $i\rhd k=k$,
$j\rhd l=l$,
\begin{align*}
\mK_Y^\veps=\ku\langle x,y,z,w\,:\, &  x^2 , y^2, z^2, w^2,  zy +\veps yz,  wx
+\veps xw, \\
& yxy -\veps xyx,  zxz -\veps xzx, wyw -\veps ywy,  wzw -\veps zwz, \\
& zxyx + yzxy,  zxyz +\veps xzxy,  \\
& wyx -\veps zwy - yxz + \veps xzw,  wzx -\veps zxy - ywz +\veps xyw, \\
&  wyzxy - \veps ywyzx - xyzwy + xyxzw, \\
& wyzxw + zxywz - yxzwy - xwyzx, \\
& wyzw -\veps zxwz - yzxw + yxwy +\veps xwyz -\veps xyzx  \rangle,
\end{align*}
and  $\dim\mK^\veps=144$.
\end{enumerate}

\noindent $\dim\mK^\veps(1)=5$:
\begin{enumerate}
\item[(9)] $Y=\{i,j,k,l,m\}$,  $i\rhd j= k$,
$i\rhd l=m$, $j\rhd l\neq l$, $k\rhd m\neq m$, $j\rhd m=m$, $k\rhd
l=l$,
\begin{align*}
\mK^\veps=\ku\langle x,y,z,w,u\, :\,  &x^2,  y^2, z^2, w^2, u^2, wz +\veps zw,
uy + \veps yu, \\
& zx +\veps yz +\veps xy,  zy + yx +\veps xz, \\
& ux +\veps wu +\veps xw,  uw + wx +\veps xu,\\
&  yxy -\veps xyx,  wxw -\veps xwx,   wyw -\veps ywy,   uzu -\veps zuz,\\
& wyx +\veps wxz -\veps zwy, wyz + wxy - zwx, \\
& uzw -\veps wxz - xuz, wxyz - zwxz,  \\
&  wxyw + ywxy + xywx, \\
& wxyxz -\veps zwxyx, wxzw + xwxz, \\
& wxyxwx +\veps ywxyxw, wxyxwy +\veps xwxyxw\rangle,
\end{align*}
and  $\dim\mK^\veps=288$.
\end{enumerate}

The stabilizers subgroups of $\s_4$ are, in each case, the following:
\begin{enumerate}
  \item $\Z_2\times \Z_2\cong  \langle  g_i , g_j\rangle\subset \s_4$  with $
i\rhd j= j$;
  \item $D_4\cong \langle  g_i, \sigma \rangle\subset \s_4$ (if, e.g.,
$g_i=(12)$, $\sigma=(1324)$);
\item $\Z_2\cong \langle  g_k\rangle\subset \s_4$, $ k=i\rhd j$.
\item $\s_3\cong\langle
 g_i , g_j, g_k\rangle\subset\s_4$, $ i \rhd j= k $;
\item
$\Z_2\cong \langle  g_j g_l\rangle$, $ j\neq l,
 j\rhd  l= l$;
\item $\s_3\cong\langle g_{i\rhd j} , g_{j\rhd k}, g_{k\rhd i}
\rangle\subset\s_4$;
\item If $\mK^\veps$ belongs to items (7) to (8) then $\stab\mK^\veps=1$.\qed
\end{enumerate}
\end{prop}

\begin{exas}
We give, as an illustration, an example of a subset $Y\subseteq \mO_2^4$ for
each case in the previous proposition. Note that any comodule algebra $\mK_{Y'}$
such that $Y'$ is not in the next list, is $\s_4$-conjugated to another,
$\mK_{Y}$, with $Y$ a set in the list.
\begin{enumerate}
\item $Y=\{(12)\}$,
\item $Y=\{(12),(34)\}$,
\item $Y=\{(12),(13)\}$,
\item $Y=\{(12),(13),(23)\}$,
\item $Y=\{(12),(13),(34)\}$,
\item $Y=\{(12),(13),(14)\}$,
\item $Y=\{(12),(13),(23),(14)\}$,
\item $Y=\{(12),(13),(24),(34)\}$,
\item $Y=\{(12),(13),(23),(14),(24)\}$.\qed
\end{enumerate}
\end{exas}

\begin{rmk}\label{rem:kYkZ}
Let $Y\subset \mO_2^4$ and let $Z\subset\mO_2^4$ be such that $\mO_2^4=Y\sqcup
Z$, as sets. Denote by $Y_j$ one of the subsets of item $j$ of Proposition
\ref{pro:coideal transp}, and by $Z_j$ the corresponding complement. Notice that
we have the following bijections
\begin{align*}
& Z_1\cong Y_9, & Z_2\cong Y_8, && Z_3\cong Y_7, && Z_4\cong Y_6, && Z_5\cong
Y_5. &
\end{align*}
Therefore, we have that $\dim\mK_Y\dim\mK_Z=\dim\B^\eps$, for every $Y$. An
analogous statement holds for the case $X=\mO_4^4$.
\end{rmk}

\section{Representations of $\Rep(\widehat{\B_2}(X,q)\# \ku
G)$}\label{main:section}

In this Section, we take $\mQ=(X,q,G,(\cdot,g,(\chi_l)_{l\in X}),(\lambda_C)_{C
\in\mR'})$ as one of ql-data from Section \ref{pointe-over-sn}. Note that in
this case, the set $ C_i=\{(i,i)\}$ belongs to $\mR=\mR'$ and $(i\rhd j)\rhd
i=j$, for any $i, j\in X$. Let $\hq$ be the corresponding Hopf algebra defined
in Section \ref{hopf-ql} and set $\mH= \widehat{\B_2}(X,q)\# \ku G$. We will
assume $\dim\widehat{\B_2}(X,q)<\infty$ (and thus $\dim\hq<\infty$,
\cite[Proposition 4.2]{GG}). In particular, this holds for $n=3,4,5$.

\subsection{$\widehat{\B_2}(X,q)\# \ku G$-Comodule algebras  }
We shall construct families of comodule algebras over quadratic approximations
of Nichols algebras. These families are large enough to classify module
categories in all of our examples.

\begin{defi}\label{def:compatible}
Let $F< G$ be a subgroup and $\psi\in Z^2(F, \ku^{\times})$. If $Y\subseteq X$
is a subset such that  $F\cdot Y\subseteq Y$, that is $F< \stab \mK_Y$, we shall
say that a family of scalars  $\xi=\{\xi_C\}_{ C\in \mR}$, $\xi_C\in \ku$ is
\emph{compatible} with the triple $(Y, F, \psi)$ if for any  $ f\in \stab
\mK_Y$,
\begin{align*}
\xi_{f\cdot C}\, \chi_i(f)\chi_j(f)&=\xi_C\, \psi(f, g_i
g_j)\,\psi(fg_ig_j,f^{-1}), && \text{ if }  C\in \mR^Y_1, (i,j)\in C;\\
\xi_{f\cdot C}\, \chi^2_i(f)\chi_j(f)&=\xi_C\,
\psi(f, g_i g_jg_i )\,\psi(fg_ig_jg_i ,f^{-1}), &&\text{ if } C\in \mR^Y_2,  (i,j)\in C;\\
\xi_{C_i}&=\xi_{C_j}=0, && \text{ if } C\in \mR^Y_2, (i,j)\in C.
\end{align*}
We will assume that the family $\xi$ is normalized by $\xi_C=0$ if either $C \in \mR^Y_1$, $(i,j)\in C$, and $g_ig_j\notin F$
or if $C \in \mR^Y_2$, $(i,j)\in C$, and $g_ig_jg_i\notin F$.
\end{defi}

We now introduce the comodule algebras we shall work with.

\begin{defi} Let  $F< G$ be a subgroup, $\psi\in Z^2(F, \ku^{\times})$,
and let $Y\subseteq X$ be a subset such that  $F\cdot Y\subseteq Y$ and let
$\xi=\{\xi_C\}_{ C\in \mR'}$ be compatible with $(Y, F, \psi)$. Define  $\Ac(Y,
F, \psi, \xi)$ to be the algebra generated by $\{y_l, e_{f} : l\in Y, f\in F\}$
and relations
\begin{align}\label{eqn:gq21} &e_1=1, \quad e_re_s=\psi(r,s)\, e_{rs}, \qquad
r,s\in F,\\
\label{eqn:gq22} &e_f\, y_l = \chi_l(f)\, y_{f\cdot l}\, e_f, \qquad f\in F, \,
l \in Y, \\
\label{eqn:gq23}&\vartheta_{C,Y}(\{y_l\}_{l\in X}) =\begin{cases} \xi_C\,
e_{C}\quad \text{ if } e_{C}\in F\\
0 \qquad\;\;\,\, \text{ if } e_{C}\notin F
 \end{cases}\quad C\in \mR.
\end{align}
Here $\vartheta_{C,Y}$ was defined in \eqref{eqn:psiC} and the element $e_C$ is
defined as
\begin{align}\label{eqn:eC}
  e_{C}=\begin{cases}
    e_{g_ig_j} &\text{if} \quad C\in \mR^Y_1, (i,j)\in C, \\
    e_{g_ig_jg_i} & \text{if} \quad C\in \mR^Y_2, (i,j)\in C \cap Y\times Y, \\
    0, & \text{if} \quad C\in \mR^Y_3.
  \end{cases}
\end{align}

If $Z\subseteq X$ is a subset invariant under the action of $F$ we define
$\Bc(Z,F,\psi,\xi)$ as the subalgebra of $\Ac(X,F,\psi,\xi)$ generated by
elements  $\{y_l, e_{f} : l\in Z, f\in F\}$.
\end{defi}

\begin{rmk} \begin{itemize}
              \item[(a)]  Applying $ad(f)$, $f\in \stab \mK_Y$ to equation
\eqref{eqn:gq23}  and using \eqref{chi-q} one can deduce the equations
in Definition \ref{def:compatible}.
              \item[(b)] It may happen that $\Bc(Z,F,\psi,\xi)\neq
\Ac(Z,F,\psi,\xi)$.
            \end{itemize}

\end{rmk}

Let $\lambda:\Ac(Y, F, \psi, \xi)\to \mH \otk
\Ac(Y, F, \psi, \xi)$ be the map defined by
\begin{align}\label{eqn:coaction}
&\lambda(e_f)=f\ot e_f,\;  \lambda(y_l)=x_l\ot 1+ g_l\ot y_l,
\end{align}
for all $f\in F$, $l\in Y$.

\begin{lema}
$\Ac(Y, F, \psi, \xi)$ is a left $\mH$-comodule algebra with coaction $\lambda$
as in \eqref{eqn:coaction} and $\Bc(Z,F,\psi,\xi)$ is a subcomodule algebra of
$\Ac(X, F, \psi, \xi)$.
\end{lema}
\pf
Let us prove first that the map $\lambda$ is well-defined. It is easy to see
that $\lambda(e_f y_l)= \chi_l(f)\, \lambda( y_{f\cdot l}\, e_g)$ for any $f\in
F, l\in X$.

Let $C\in  \mR^Y_1$ and $(i,j)\in C$. In this case $\vartheta_C=\phi_C$. We
shall prove that
$\lambda(\phi_C(\{y_l\}_{l\in X}))=\lambda(\xi_C\, e_{g_ig_j})$.
Using the definition of the polynomial $\phi_C$ we obtain that
\begin{align*} \lambda(\phi_C(\{y_l\}_{l\in X}))&=\sum_{h=1}^{n(C)}\eta_h(C) \,
x_{i_{h+1}} x_{i_{h}}\ot 1 +
x_{i_{h+1}} g_{i_h}\ot y_{i_h}
+ g_{i_{h+1}} x_{i_h}\ot y_{i_{h+1}}+\\
& + g_{i_{h+1}}g_{i_{h}}\ot y_{i_{h+1}} y_{i_{h}}=\phi_C(\{x_l\}_{l\in X})\ot 1
+ g_i g_j \ot \phi_C(\{y_l\}_{l\in X})\\
&= \xi_C\; g_i g_j \ot e_{g_ig_j}=\lambda(\xi_C\, e_{g_ig_j}).
\end{align*}

The second equality follows since $i_{n(C)+1}=i_1$ and
$$ g_{i_{h+1}} x_{i_h}= q_{i_{h+1}i_{h}}\; x_{i_{h+2}}g_{i_{h+1}}, \quad
\eta_h(C) q_{i_{h+1}i_{h}}=- \eta_{h+1}(C).$$
\smallbreak

Now, let $C\in  \mR^Y_2$, $(i,j)\in C $ and   $i\rhd j\notin Y$. In this case
relation \eqref{eqn:gq23} is $y_iy_jy_i+q_{i\rhd j i}y_jy_iy_j= \xi_C\, e_{g_ig_jg_i}$. Note
that assumption $\xi_{C_i}=\xi_{C_j}=0$ implies that $y^2_i=0=y^2_j$. The proof
that $\lambda(y_iy_jy_i+q_{i\rhd j i}y_jy_iy_j)=\xi_C\lambda(e_{g_ig_jg_i})$ is a
straightforward computation.
\epf

\begin{teo}\label{comod-alg1}
Let $Y\subseteq X$ be an $F$-invariant subset and assume that $\Ac(X, F, \psi,
\xi)\neq 0$, then the following statements hold:
\begin{itemize}
 \item[1.] The algebras $\Ac(X, G, \psi, \xi)$ are left $\mH$-Galois extensions.
 \item[2.] If $\xi$ satisfies
 \begin{equation}\label{scalars-bigalois} \xi_C=\begin{cases} - \lambda_C &
\text{ if }  \quad  \lambda_C\neq 0, \\
  \quad 0   & \text{ if }  \quad \lambda_C=0 \text{ and } g_j g_i\neq 1,\\
  \text{ arbitrary } & \text{ if }  \quad \lambda_C=0 \text{ and } g_j g_i= 1.
  \end{cases}
  \end{equation}
 then $\Ac(X, G, 1, \xi)$ is a  $(\mH,\hq)$-biGalois object.

   \item[3.] $\Bc(Y,F,\psi,\xi)_0=\ku_{\psi} F$ and thus $\Bc(Y,F,\psi,\xi)$ is
a right $\mH$-simple left $\mH$-comodule algebra.
 \item[4.] There is an isomorphism of comodule algebras $\gr \Bc(Y, F, \psi,
\xi)\simeq \mK_Y\# \ku_{\psi} F$.
\item[5.]  There is an isomorphism $\Bc(Y, F, \psi, \xi)\simeq \Bc(Y', F',
\psi', \xi')$ of comodule algebras if and only if
      $Y=Y'$, $F=F'$, $\psi=\psi'$ and $\xi=\xi'$.

\end{itemize}
\end{teo}

\pf 1.  To prove that  $\Ac(X, G, \psi, \xi)$ is a Galois observe
 that the canonical map $$\can:\Ac(X, G, \psi, \xi)\otk \, \Ac(X, G, \psi,
\xi)\to \mH\, \otk \, \Ac(X, G, \psi, \xi),$$ $\can(x\ot y)=x\_{-1}\ot x\_0 y$,
is surjective. Indeed for any $f\in G$, $l\in X$
$\can(e_f\ot e_{f^{-1}})=f\ot 1,$  $\can(y_l\ot 1- e_{g_{l}}\ot e_{g_l^{-1}}
y_l)=x_l\ot 1$.

\medbreak

2. Define the map
$
\rho:\Ac(X, G, 1, \xi)\to \Ac(X, G, 1, \xi)\otk  \hq,$ by
\begin{align*}
&\rho(e_f)=e_f\ot H_f,\; \rho(y_l)=y_l\ot 1+ e_{g_l}\ot a_l,\quad l\in X, f\in
G.
\end{align*}
The map $\rho$ is well-defined. Indeed, if $C\in  \mR$ and $(i,j)\in C$ then
\begin{align*} \rho(\phi_C(\{y_l\}_{l\in X}))&= \phi_C(\{y_l\}_{l\in X})\ot 1 +
e_{g_ig_j}\ot \phi_C(\{a_l\}_{l\in X})\\
&= \xi_C\, e_{g_ig_j} \ot 1 +  \lambda_C\,
e_{g_ig_j}\ot \big( 1- H_{g_ig_j}\big).
\end{align*}
Clearly if $\xi$ satisfies \eqref{scalars-bigalois} then
$\rho(\phi_C(\{y_l\}_{l\in X}))= \xi_C\,\rho( e_{g_ig_j})$. The
proof that $\Ac(X, G, 1, \xi)$ is a $(\mH,\hq)$-bicomodule and a right
$\hq$-Galois object is done by a straightforward computation.

\medbreak

3. If $\Ac(X, F, \psi, \xi)\neq 0$ then there is a group $\overline{F}$ with a
projection $F\twoheadrightarrow \overline{F}$ such that $\Ac(Y, F, \psi,
\xi)_0=\ku_{\overline{\psi}} \overline{F}.$ The map $\Ac(Y, F, \psi, \xi)_0\ot
\Ac(Y, F, \psi, \xi)_0 \to \ku F\ot \Ac(Y, F, \psi, \xi)_0$ defined by $e_f\ot
e_g\mapsto f\ot \psi(f,g)\, e_{fg}$ is surjective. Hence $F=\overline{F}$. This
implies that $\Bc(Z,F,\psi,\xi)_0=\ku_{\psi} F$ and by \cite[Prop. 4.4]{M1}
follows that $\Bc(Z,F,\psi,\xi)$ is a right $\mH$-simple left $\mH$-comodule
algebra.

 \medbreak
 4.
Follows from Theorem \ref{mod-over-pointed} (3) that $\gr \Bc(Y, F, \psi,
\xi)\simeq \mK\# \ku_{\psi} F$ for some homogeneous left coideal subalgebra
$\mK\subseteq \widehat{\B_2}(X,q)$. Recall that $\mK$ is identified with the
subalgebra of $\gr\Bc(Y, F, \psi, \xi)$ given by $$\{a\in \gr\Ac(Y, F, \psi,
\xi) : (\id\otimes \pi)\lambda(a)\in \mH\ot 1\}.$$ See \cite[Proposition 7.3
(3)]{M1}. In \textit{loc. cit.} it is also proved that the composition
$$\gr\Bc(Y, F, \psi, \xi) \xrightarrow{(\theta\ot \pi)\lambda}  \mK\# \ku_{\psi}
F \xrightarrow{\;\mu\;} \gr\Bc(Y, F, \psi, \xi), $$
is the identity map, where $\theta:\mH\to \widehat{\B_2}(X,q)$, $\pi:\gr\Bc(Y,
F, \psi, \xi)\to \ku_{\psi} F$ are the canonical projections and $\mu$ is the
multiplication map. Both maps are bijections and since for any $l\in Y$, $
(\theta\ot \pi)\lambda(y_l)=x_l$, then $\mK=\mK_{Y}$.

 \medbreak

 5. Let $\beta:\Bc(Y, F, \psi, \xi)\to \Bc(Y', F', \psi', \xi')$ be
a comodule algebra isomorphism. The restriction of $\beta$ to $\Bc(Y, F, \psi,
\xi)_0$ induces an isomorphism between $\ku_\psi F$ and $\ku_{\psi'} F'$, thus
$F=F'$ and $\psi=\psi'$. Since $\beta$ is a comodule morphism it is clear that
$Y=Y'$ and $\xi_C=\xi'_C$ for any $C\in\mR$.
\epf

\begin{cor}\label{cocycle-def} If $\Ac(X, G, 1, \xi)\neq 0$ for some $\xi$
satisfying \eqref{scalars-bigalois}, then
\begin{enumerate}
  \item[1.] The Hopf algebras $\mH=\widehat{\B_2}(X,q)\# \ku G $ and $\hq$ are
cocycle deformations of each other.
  \item[2.]  There is a bijective correspondence between equivalence classes of
exact module categories over $\Rep(\mH)$ and $\Rep(\hq)$.\qed
\end{enumerate}
\end{cor}

\begin{rmk}
Under the assumptions in Corollary \ref{cor:cocycle def}, we obtain,
in particular, that $\gr\hq=\widehat{\B_2}(X,q)\#\ku G$, since the latter is a
quotient of the first.
\end{rmk}

The following corollary uses Propositions \ref{pro:3nonull} and \ref{pro:4nonull},
where certain algebras are shown to be not null. These propositions will be proven in the
Appendix and their proofs are independent of the results in the article.

\begin{cor}\label{cor:cocycle def}
Let $H$ be a non-trivial pointed Hopf algebra over $\s_3$ or $\s_4$. Then $H$ is
a cocycle deformation of $\gr H$.
\end{cor}
\begin{proof}
Finite-dimensional Nichols algebras over $\s_3$ and $\s_4$ coincide with their
quadratic approximations. That is, if $H$ is a finite-dimensional pointed Hopf
algebra over $\s_n$, $n=3,4$, then $\gr H\cong \widehat{\B_2}(X,q)\#\ku\s_n$. By
\cite[Main Theorem]{GG} we know that $H\cong\hq$. Therefore, the theorem follows
from Corollary \ref{cocycle-def}, since in Propositions \ref{pro:3nonull},
\ref{pro:4nonull} we show the existence of non-zero $(\gr\hq,\hq)$-biGalois objects in these cases.

When dealing with $\mQ_4^{-1}[t]$ or $\mD[t]$, notice that condition
$\xi_2=2\xi_1$ in Proposition \ref{pro:4nonull} does not interfere with the
proof, since, by equation \eqref{scalars-bigalois}, $\xi_1$, resp. $\xi_2$, can
be chosen arbitrarily.
\end{proof}

\begin{rmk} In \cite[Theorem A1]{Ma} Masuoka proved that the Hopf algebras
$u(\Do, \lambda,\mu)$ associated to a datum of finite Cartan
type $\Do$  appearing in the classification of Andruskiewitsch and Schneider
\cite{AS} are cocycle deformations to the associated graded Hopf algebras
$u(\Do, 0,0)$.

Corollaries \ref{cocycle-def} (1) and \ref{cor:cocycle def} provide a similar
result for some families of Hopf algebras constructed from  Nichols algebras not
of diagonal type. It would be interesting to generalize this kind of result for
larger classes of Nichols algebras.
\end{rmk}

\subsection{Module categories over $\Rep(\hq)$}

\

Let $A$ be a $\mH$-comodule algebra with $\gr A=\mK_Y\#\ku_\psi F$, for
$F\leq \stab\mK_Y$, $\psi\in Z^2(F,\ku^*)$. Let $Z$ be such that $X=Y\sqcup Z$
as sets. Notice that $F\leq\stab\mK_Z$.

\begin{lema}\label{lem:subalgebra}
Under the above assumptions there exists a family of scalars $\xi$ compatible
with $(X,F,\psi)$ such that $A\simeq \Bc(Y,F,\psi,\xi)$ as comodule algebras.
\end{lema}
\pf
The canonical projection $\pi:A_1\to A_1/A_0\simeq \mK_Y(1)=\ku Y$ is a morphism
of $\Ac_0$-bimodules. Let $\iota:\ku Y\to  A_1$ be a section of
$\Ac_0$-bimodules of $\pi$. Since elements $\{x_l: l\in Y\}$ are in the image of
$\pi$ we can choose elements $\{y_l: l\in Y\}$ in $A_1$ such that
$\iota(x_l)=y_l$ for any $l\in Y$. It is straightforward to verify that
\begin{align*} \lambda(y_l)=x_l\ot 1+ g_l\ot y_l, \quad e_f\, y_l = \chi_l(f)\,
y_{f\cdot l}\, e_f, \qquad f\in F, \, l \in Y.
\end{align*}
Since $\gr A$ is generated by elements $\{x_l, e_{f} : l\in Y, f\in F\}$ then
$A$ is generated as an algebra by elements $\{y_l, e_{f} : l\in Y, f\in F\}$.

Now, let $B=A\ot \mK_Z$. Then $B$ has an comodule algebra structure for which
the canonical inclusion $A\hookrightarrow A\ot 1\subset B$ is a homomorphism.
The algebra structure is given as follows. For $i\in Y$, $j\in Z$, $f\in F$,
\begin{align*}
(y_i\ot 1)(1\ot y_j)&=(y_i\ot y_j);\\
(1\ot y_j)(y_i\ot 1)&=\begin{cases}
                      q_{ji}y_i\ot y_j+\xi_Ce_C\ot 1, & \text{if } i\rhd j=j\\
q_{ji}y_{j\rhd i}\ot y_j-q_{ji}q_{j\rhd i\,j}y_iy_{j\rhd i}\ot 1+\xi_Ce_C\ot 1,
& \text{if } i\rhd j\neq j, \\
& \quad i\rhd j\in Y;\\
q_{ji}1\ot y_{j\rhd i}y_j-q_{ji}q_{j\rhd i\,j}y_i\ot y_{j\rhd i}+\xi_Ce_C\ot 1,
& \text{if } i\rhd j\neq j,\\
& \quad i\rhd j\notin Y;
                      \end{cases}\\
(e_f\ot 1)(1\ot y_j)&=e_f\ot y_j;\\
(1\ot y_j)(e_f\ot 1)&=\chi_j^{-1}(f)e_f\ot y_{f^{-1}\cdot j}.
\end{align*}
Here $C$ stands for the class $C\in\mR'$ such that $(j,i)\in C$. Recall that by
definition $\xi_C=0$ if $g_C\notin F$. Then the map
\begin{equation}\label{eqn:multiplication}
m:B\to \Ac(X,F,\psi,\xi), \quad a\ot x\mapsto ax
\end{equation}
is an algebra epimorphism. Now, if
\begin{align*}
A\ni a&\mapsto a_{(-1)}\ot a_{(0)}\in \mH\ot
A \quad \text{ and }\\
\mK_Z\ni x&\mapsto x_{(-1)}\ot x_{(0)}\in \mH\ot\mK_Z
\end{align*}
denote the
corresponding coactions, define $\lambda:B\to \mH\ot B$ by $\lambda(a\ot
x)=a_{(-1)}x_{(-1)}\ot a_{(0)}\ot x_{(0)}$. It is straightforward to check that
$\lambda$ is well defined. We do this case by case in the definition of the
multiplication of $B$ above. For instance, if $i\rhd j\neq j$ and $i\rhd j\in
Y$, then we have
\begin{align*}
\lambda(1\ot y_j)\lambda(y_i\ot 1)&=(g_j\ot(1\ot y_j)+x_j\ot (1\ot
1))(g_i\ot(y_i\ot 1)+x_i\ot (1\ot 1))\\
&=(g_j\ot(1\ot y_j))(g_i\ot(y_i\ot 1))+(x_j\ot (1\ot 1))(g_i\ot(y_i\ot 1))\\
&\quad+(g_j\ot(1\ot y_j))(x_i\ot (1\ot 1))+(x_j\ot (1\ot 1))(x_i\ot (1\ot 1))\\
&=g_jg_i\ot (1\ot y_j)(y_i\ot 1)+x_jg_i\ot (y_i\ot 1)\\
&\quad +q_{ji}x_{j\rhd i}g_j\ot(1\ot y_j)+x_jx_i\ot (1\ot 1)\\
&=g_jg_i\ot (q_{ji}y_{j\rhd i}\ot y_j-q_{ji}q_{j\rhd i\,j}y_iy_{j\rhd i}\ot
1+\xi_Cg_C\ot 1)\\
&\quad +x_jg_i\ot (y_i\ot 1)+q_{ji}x_{j\rhd i}g_j\ot(1\ot y_j)\\
&\quad +(q_{ji}x_{j\rhd i}x_j-q_{ji}q_{j\rhd i\,j}x_ix_{j\rhd i}\ot 1)\ot (1\ot
1),
\end{align*}
which coincides with $\lambda(q_{ji}y_{j\rhd i}\ot y_j-q_{ji}q_{j\rhd
i\,j}y_iy_{j\rhd i}\ot 1+\xi_Cg_C\ot 1)$.

Thus, $B$ is an $\mH$-comodule algebra with $$\dim B=\dim
A\dim\mK_Z=\dim\mK_Y\dim\mK_Z|F|=\dim\Ac(X,F,\psi,\xi),$$ by Remark
\ref{rem:kYkZ} and then the map $m$ from \eqref{eqn:multiplication} is an
isomorphism.
\epf

We can now formulate the main result of the paper. For any $h\in G$, we denote
$\xi^h_C=\xi_{h^{-1}\cdot C}$. Recall that we denote by $\Bc(Y, F, \psi, \xi)$
the sub-comodule algebra of $\Ac(X, F, \psi, \xi)$ generated by $\{y_i\}_{i\in
Y}$.

\begin{teo}\label{clasification-general}
\begin{itemize}
\item[1.] Let $\Mo$ be an exact indecomposable module category over $\Rep(\hq)$,
then there exists
\begin{itemize}
  \item[(i)] a subgroup $F< G$, and a 2-cocycle $\psi\in Z^2(F,\ku^{\times})$,
  \item[(ii)] a subset $Y\subset X$ such that $F\cdot Y\subset Y$,
  \item[(iii)] a family of scalars $\{\xi_C\}_{ C\in \mR'}$ compatible with $(X,
F,\psi)$,
\end{itemize}
such that there is a module equivalence $\Mo\simeq {}_{\Bc(Y, F, \psi,
\xi)}\Mo$.

\smallbreak

\item[2.] Let $(Y, F, \psi, \xi)$,
$(Y', F', \psi', \xi')$ be two families as before. Then there is an equivalence
of module categories ${}_{\Bc(Y, F, \psi, \xi)}\Mo\simeq {}_{\Bc(Y', F', \psi',
\xi')}\Mo$ if and only if there exists an element $h\in G$ such that $F'=h F
h^{-1}$, $\psi'=\psi^h$, $Y'=h\cdot Y$ and $\xi'=\xi^h$.\end{itemize}
\end{teo}

\pf 1.  By Corollary \ref{cor:cocycle def} we can assume that $\Mo$ is an exact
indecomposable module category over $\gr\hq=\mH$.  It follows by \cite[Theorem
3.3]{AM} that there is a right $\mH$-simple left $\mH$-comodule algebra $\Ac$
such that $\Mo\simeq {}_{\Ac}\Mo$.  Theorem \ref{mod-over-pointed} implies that
there is a subgroup $F< G$, a 2-cocycle $\psi\in Z^2(F,\ku^{\times})$ and  a
subset $Y\subset X$ with $F\cdot Y\subset Y$ such that $\gr\Ac= \mK_Y\# \ku_{\psi}
F$. Here $\Ac_0= \ku_{\psi} F$. Then the result follows from Lemma
\ref{lem:subalgebra}.

\medbreak

2. Assume that the module categories  ${}_{\Bc(Y, F, \psi, \xi)}\Mo, {}_{\Bc(Y',
F', \psi', \xi')}\Mo$ are equivalent, then  Theorem \ref{Morita-equivalence}
implies that there exists an element $h\in G$ such that $\Bc(Y', F', \psi',
\xi')\simeq h\Bc(Y, F, \psi, \xi)h^{-1}$ as $H$-comodule algebras.

The algebra map $\alpha: h\Bc(Y, F, \psi, \xi)h^{-1}\to \Bc(h\cdot Y, hFh^{-1},
\psi^h, \xi^h)$ defined by
$$\alpha(he_fh^{-1})= e_{hfh^{-1}},\quad  \alpha(hy_lh^{-1})= \chi_l(h)
\,y_{h\cdot l},$$
for any $f\in F$, $l\in Y$, is a well-defined comodule algebra isomorphism.
Whence   $\Bc(Y', F', \psi', \xi')\simeq \Bc(h\cdot Y, hFh^{-1}, \psi^h, \xi^h)$
and using Theorem \ref{comod-alg1} (3) we get the result.\epf

As a consequence of Theorem \ref{clasification-general} we have the following
result.
\begin{cor}\label{hopfgalois}  Any $\mH$-Galois object is of the form $\Ac(X, G,
\psi, \xi)$.
\end{cor}
\pf
Let $A$ be a $\mH$-Galois object. Then ${}_A\Mo$ is an exact module category
over $\Rep \mH$. Moreover, ${}_A\Mo$ is indecomposable. In fact, otherwise there
would exist a proper bilateral ideal $J\subset A$ $\mH$-stable \cite[Proposition
1.18]{AM}. Thus, $\can(A\ot J)=\can(J\ot A)$, what contradicts the bijectivity
of $\can$. Then, by Theorem \ref{clasification-general} there exists a datum
$(X, G, \psi, \xi)$ such that $A\cong \Ac(X, G, \psi, \xi)$.
\epf

\subsection{Modules categories over $\B(\mO_2^3,-1)\#\ku\s_3$}

We apply Theorem \ref{clasification-general} to exhibit explicitly all module
categories in this particular case. In this case the rack is
$$\mO_2^3=\{(12),(13),(23)\}.$$
For each $i\in \mO_2^3$ we shall denote by $g_i$ the  element $i$ thought as an
element in the group $\s_3$. We will show in the Appendix that the algebras in
the following result are not null. Then the next corollary follows from Theorem
\ref{clasification-general}.

\begin{cor}\label{clasif:paricular} Let $\Mo$ be an indecomposable exact module
category over $\Rep(\B(\mO_2^3,-1)\#\ku\s_3)$. Then there is a module
equivalence $\Mo\simeq {}_\Ac\Mo$ where $\Ac$ is one (and only one) of the
comodule algebras in following list. In the following $i,j, k$ denote elements
in $\mO_2^3$ and $\xi, \mu, \eta\in \ku$.

\begin{enumerate}

\item[1.] For any subgroup $F\subseteq \s_3$, $\psi\in Z^2(F,\ku^\times)$, the
twisted group algebra $\ku_\psi F$.

  \item[2.] The algebra $\Ac(\{i\},\xi, 1)=< y_i: y_i^2=\xi 1 > $, with coaction
determined by $\lambda(y_i)= x_i\ot 1+ g_i\ot y_i$.

  \item[3.] The algebra  $\Ac(\{i\},\xi, \Z_2)=< y_i, h: y_i^2=\xi 1, h^2=1, h
y_i=-y_i h > $ with coaction determined by $\lambda(y_i)= x_i\ot 1+ g_i\ot y_i$,
$\lambda(h)=g_i\ot h$.

     \item[4.] The algebra  $\Ac(\{i,j\}, 1)=< y_i, y_j: y_i^2=y_j^2=0,\;
y_iy_jy_i=y_jy_iy_j > $ with coaction determined by $\lambda(y_i)= x_i\ot 1+
g_i\ot y_i$, $\lambda(y_j)= x_j\ot 1+ g_j\ot y_j.$

\item[5.]  The algebra $\Ac(\{i,j\}, \Z_2)=< y_i, y_j,h: y_i^2=y_j^2=0, h^2=1,
hy_i=-y_j h, y_iy_jy_i=y_jy_iy_j > $ with coaction determined by $\lambda(y_i)=
x_i\ot 1+ g_i\ot y_i$, $\lambda(y_j)= x_j\ot 1+ g_j\ot y_j,$ $\lambda(h)=g_k\ot
h$, where $k\neq i,j$.

\item[6.] The algebra $ \Ac(\mO_2^3,\xi, 1)$, generated by $\{y_{(12)},
y_{(13)},$ $ y_{(23)}\}$ subject to relations
\begin{align*}
&y_{(12)}^2=y_{(13)}^2=y_{(23)}^2=\xi 1,\\
&y_{(12)}y_{(13)}+y_{(13)}y_{(23)}+y_{(23)}y_{(12)}=0,\\
&y_{(13)}y_{(12)}+y_{(23)}y_{(13)}+y_{(12)}y_{(23)}=0.
\end{align*}
The coaction is determined by $\lambda(y_s)= x_s\ot 1+ g_s\ot y_s$ for any $s\in
\mO_2^3$.

\item[7.] The algebra $ \Ac(\mO_2^3,\xi, \Z_2)$, generated by $\{y_{(12)},
y_{(13)}, y_{(23)},$ $ h\}$ subject to relations
\begin{align*}
&y_{(12)}^2=y_{(13)}^2=y_{(23)}^2=\xi 1,\;\; h^2=1,\\
&h y_{(12)}= -y_{(12)} h,\;\; h y_{(13)}= -y_{(23)} h,\\
&y_{(12)}y_{(13)}+y_{(13)}y_{(23)}+y_{(23)}y_{(12)}=0.
\end{align*}
The coaction is determined by $\lambda(h)=g_{(12)}\ot h$, $\lambda(y_s)= x_s\ot
1+ g_s\ot y_s$ for any $s\in \mO_2^3$.

\item[8.] The algebra $ \Ac(\mO_2^3,\xi, \mu, \eta, \Z_3)$, generated by $\{
y_{(12)}, y_{(13)}, $ $y_{(23)}, h\}$ subject to relations
\begin{align*}
& y_{(12)}^2=y_{(13)}^2=y_{(23)}^2=\xi 1,\;\; h^3=1,\\
& h y_{(12)}= y_{(13)} h,\;\;  h y_{(13)}= y_{(23)} h,\;\;  h y_{(23)}= y_{(12)}
h,\\
&y_{(12)}y_{(13)}+y_{(13)}y_{(23)}+y_{(23)}y_{(12)}=\mu\, h,\\
&y_{(13)}y_{(12)}+y_{(23)}y_{(13)}+y_{(12)}y_{(23)}=\eta\, h^2.
\end{align*}
The coaction is determined by $\lambda(h)=g_{(132)}\ot h$, $\lambda(y_s)= x_s\ot
1+ g_s\ot y_s$, for any $s\in \mO_2^3$.

\item[9.] The algebras $ \Ac(\mO_2^3,\xi, \mu, \s_3, \psi)$, for each  $\psi\in
Z^2(\s_3,\ku^\times)$, generated by $\{y_{(12)}, y_{(13)}, y_{(23)},
e_h\,:\,h\in\s_3\}$  subject to relations
\begin{align*}
&e_h e_t= \psi(h,t)\, e_{ht},\; \;   e_h y_s=- y_{h\cdot s} e_h \quad\; \,
h,t\in \s_3, s\in \mO_2^3,\\
& y_{(12)}^2=y_{(13)}^2=y_{(23)}^2=\xi 1,\\
&y_{(12)}y_{(13)}+y_{(13)}y_{(23)}+y_{(23)}y_{(12)}=\mu\, e_{(123)}.
\end{align*}
The coaction is determined by $\lambda(e_h)=h\ot e_h$, $\lambda(y_s)= x_s\ot 1+
g_s\ot y_s$ for any $s\in \mO_2^3$.\hfill \qed
\end{enumerate}
\end{cor}

\section{Appendix: $\mA(Y,F,\psi,\xi)\neq 0$}

In this part, we will complete the proofs of Corollaries \ref{cor:cocycle def}
and \ref{clasif:paricular}, by showing that the algebras involved in their
statements are not null.

\begin{prop}\label{pro:3nonull}
Let $\mA(Y,F,\psi,\xi)$ be one of the algebras in Corollary
\ref{clasif:paricular}. Then $\mA(Y,F,\psi,\xi)\neq 0$.
\end{prop}
\pf
The case $Y\neq \mO_2^3$ is clear. Set $Y=\mO_2^3$. Note that each one of these
algebras is naturally a right $\ku F$-module via $a\leftharpoonup t=ae_t$,
$a\in\mA(Y,F,\psi,\xi)$, $t\in F$. Thus, we can consider the induced
representation $W=\mA(Y,F,\psi,\xi)\ot_{\ku F}W_\eps$, where $W_\eps=\ku\{ z\}$
is the trivial $\ku F$-module. Let
\begin{align*}
B=\{1,&y_{(12)},y_{(13)},y_{(23)},y_{(13)}y_{(12)},y_{(12)}y_{(13)},y_{(12)}y_{
(23)},y_{(13)}y_{(23)},\\
&y_{(12)}y_{(13)}y_{(23)},y_{(13)}y_{(12)}y_{(23)},y_{(12)}y_{(13)}y_{(12)},
y_{(12)}y_{(13)}y_{(12)}y_{(23)}\}
\end{align*}
and consider the linear subspace $V$ of $W$ generated by $B\ot z$. We show that
this is a non-trivial submodule in the four cases left, namely $F=1, \Z_2, \Z_3$
or $\s_3$. In all of the cases, the action of $y_{(12)}$ is determined by the
matrix
\begin{align*}
y_{(12)}=\left[\begin{smallmatrix}
  0& \xi& 0& 0& 0& 0& 0& 0& 0& 0& 0& 0\\1& 0& 0& 0& 0& 0& 0& 0& 0& \
0& 0& 0\\0& 0& 0& 0& 0& \xi& 0& 0& 0& 0& 0& 0\\0& 0& 0& 0& 0& \
0& \xi& 0& 0& 0& 0& 0\\0& 0& 0& 0& 0& 0& 0& 0& 0& 0& \xi& 0 \\
0& 0& 1& 0& 0& 0& 0& 0& 0& 0& 0& 0\\0& 0& 0& 1& 0& 0& 0& 0& 0& 0& \
0& 0\\0& 0& 0& 0& 0& 0& 0& 0& \xi& 0& 0& 0\\0& 0& 0& 0& 0& 0& \
0& 1& 0& 0& 0& 0\\0& 0& 0& 0& 0& 0& 0& 0& 0& 0& 0& \xi\\0& 0& \
0& 0& 1& 0& 0& 0& 0& 0& 0& 0\\0& 0& 0& 0& 0& 0& 0& 0& 0& 1& 0& 0
      \end{smallmatrix}\right].
\end{align*}

Now, take $F=\s_3$, $\psi\equiv1$. The action of $e_{(12)}$ and $e_{(13)}$ is
determined, respectively, by the matrices:
\begin{align*}
\left[\begin{smallmatrix}
  1& 0& 0& 0& \mu& 0& 0& \mu& 0& 0& 0& 0\\0& -1& 0& 0& 0& 0& 0&
0& -\mu& 0& -\mu& 0\\0& 0& 0& -1& 0& 0& 0& 0& 0& \mu& \xi&
0\\0& 0& -1& 0& 0& 0& 0& 0& \xi& -\mu& 0& 0\\0& 0& 0& 0& 0&
0& 0& -1& 0& 0& 0& 0\\0& 0& 0& 0& -1& 0& 1& 0& 0& 0& 0& -\mu\\
0& 0& 0& 0& 0& 1& 0& -1& 0& 0& 0& \mu\\0& 0& 0& 0& -1& 0& 0& 0&
0& 0& 0& 0\\0& 0& 0& 0& 0& 0& 0& 0& 0& 0& 1& 0\\0& 0& 0& 0& 0& 0&
0& 0& 0& -1& 0& 0\\0& 0& 0& 0& 0& 0& 0& 0& 1& 0& 0& 0\\0& 0& 0&
0& 0& 0& 0& 0& 0& 0& 0& 1
      \end{smallmatrix}\right] \text{and}
\left[\begin{smallmatrix}
  1& 0& 0& 0& 0& \mu& \mu& 0& 0& 0& 0& 0 \\0& 0& 0& -1& 0& 0& 0&
0& \mu& 0& \xi& 0\\0& 0& -1& 0& 0& 0& 0& 0& 0& -\mu& -\mu&
0\\0& -1& 0& 0& 0& 0& 0& 0& -\mu& \xi& 0& 0\\0& 0& 0& 0& 0&
-1& 0& 1& 0& 0& 0& -\mu\\0& 0& 0& 0& 0& 0& -1& 0& 0& 0& 0& 0\\
0& 0& 0& 0& 0& -1& 0& 0& 0& 0& 0& 0\\0& 0& 0& 0& 1& 0& -1& 0& 0&
0& 0& \mu\\0& 0& 0& 0& 0& 0& 0& 0& -1& 0& 0& 0\\0& 0& 0& 0& 0&
0& 0& 0& 0& 0& 1& 0\\0& 0& 0& 0& 0& 0& 0& 0& 0& 1& 0& 0\\0& 0& 0&
0& 0& 0& 0& 0& 0& 0& 0& 1
      \end{smallmatrix}\right].
\end{align*}
The action of $e_{(23)}$ is given by $e_{(12)}e_{(13)}e_{(12)}$. Finally, we use computer program \texttt{Mathematica}\textcopyright \ to check
that these matrices satisfy the relations defining the algebra on each case.

We deal now with a generic 2-cocycle $\psi\in Z^2(\s_3,\ku^\times)$. Let us fix $\mA=\mA(Y,F,1,\xi)$, $\mA'=\mA(Y,F,\psi,\xi)$. Also, set $U=\mK_Y\#\ku F$, $U'=\mK_Y\#\ku_\psi F$. If $\overline\psi\in Z^2(U)$ is the 2-cocycle such that $\overline\psi_{F\times F}=\psi$, see Lemma \ref{lem:extension cociclo}, it follows that $U'=U^{\overline\psi}$. Now, as $\mA$ is an $U$-comodule algebra, isomorphic to $U$ as $U$-comodules, it follows that there exists a 2-cocycle $\gamma\in Z^2(U)$ such that $\mA\cong\,_{\gamma}U$, see \cite[Sections 7 \& 8]{Mo}. It is easy to check then that $\mA'=\,_{\gamma}U'$, by computing the multiplication on the generators, and thus $\mA'\neq 0$.
\epf

To finish the proof of Corollary \ref{cor:cocycle def}, we present three
families of non trivial algebras $\Ac(X, G, 1, \xi)$, for $X=\mO_2^4$, $G=\s_4$
and certain collections of scalars $\{\xi_C\}_{C\in\mR'}$ satisfying
\eqref{scalars-bigalois}. We show $\Ac(X, G, 1, \xi)\neq 0$ in Proposition
\ref{pro:4nonull}.

\begin{defi}\label{def:extensiones s4} Let $\psi\in Z^2(\s_4,\ku^\times)$,
$\alpha,\beta\in\ku$.
\begin{itemize}
\item[1.] $\mA^{-1}_\psi(\alpha,\beta)$ is the algebra generated by $\{y_i,
e_{g} : i\in \mO_2^4, g\in \s_4\}$ with relations
\begin{align*} \quad\quad&e_1=1, \quad e_re_s=\psi(r,s)\, e_{rs}, \quad r,s\in
\s_4,\\
&e_g\, y_l = \sg(g)\, y_{g\cdot l}\, e_g, \quad\quad\quad\quad g\in \s_4, \, l
\in \mO_2^4, \\
&y_{(12)}^2= \alpha\, 1,\\
& y_{(12)} y_{(34)}+y_{(34)} y_{(12)}= 2\alpha\, e_{(12)(34)}, \\
&y_{(12)} y_{(23)}+y_{(23)} y_{(13)}+y_{(13)}y_{(12)} = \beta\, e_{(132)}.
\end{align*}

\item[2.] $\mA_\psi^4(\alpha,\beta)$ is the algebra generated by $\{y_i,
e_{g} : i\in \mO_4^4, g\in \s_4\}$ with relations
\begin{align*}
\quad\quad&e_1=1, \quad e_re_s=\psi(r,s)\, e_{rs}, \quad  r,s\in \s_4,\\
&e_g\, y_l = \sg(g)\, y_{g\cdot l}\, e_g,  \quad\quad\quad\quad g\in \s_4, \, l
\in \mO_4^4, \\
&y_{(1234)}^2= \alpha\, e_{(13)(24)},\\
& y_{(1234)} y_{(1432)}+y_{(1432)} y_{(1234)}= 2\alpha\, 1,\\
&y_{(1234)} y_{(1243)}+y_{(1243)} y_{(1423)}+y_{(1423)}y_{(1234)} =\beta\,
e_{(132)}.
\end{align*}
\item[3.] $\mA^\chi_\psi(\alpha,\beta)$ is the algebra generated by $\{y_i,
e_{g} : i\in \mO_2^4, g\in \s_4\}$  with relations
\begin{align*}
\quad\quad&e_1=1, \quad e_re_s=\psi(r,s)\, e_{rs}, \quad  r,s\in \s_4,\\
&e_g\, y_l = \chi_l(g)\, y_{g\cdot l}\, e_g,  \quad\quad\quad\quad\quad g\in
\s_4, \, l \in \mO_2^4, \\
&y_{(12)}^2= \alpha\, 1, \\
& y_{(12)} y_{(34)}-y_{(34)} y_{(12)}=0,
 \\
&y_{(12)} y_{(23)}-y_{(23)} y_{(13)}-y_{(13)}
y_{(12)} = \beta\, e_{(132)}.
\end{align*}
\end{itemize}
\end{defi}

\begin{rmk}\label{rem:extensiones s4}
Let $\mQ=\mQ^{-1}[t]$. It is clear
$\mA^{-1}_\psi(\alpha,\beta)\cong\Ac(\mO_2^4,\s_4,\psi,\xi)$ for the family
$\xi=\{\xi_C\}_{C\in\mR}$ where $\xi_C=\xi_i$, if $i=1,2,3$, is constant in the
classes $C$ with the same cardinality $|C|=i$ and where in this case
$\xi_1=\alpha$, $\xi_2=2\alpha$, $\xi_3=\beta$.

Analogously, if $\mQ=\mQ^{\chi}[t]$, $\mA^{\chi}_\psi(\alpha,\beta)$ is the
algebra $\Ac(\mO_2^4,\s_4,\psi,\xi)$ for certain family  $\xi$ subject to
similar conditions as in the previous paragraph. The same holds for
$\mQ=\mD[t]$, $\mA^4_\psi(\alpha,\beta)$ and  $\Ac(\mO_4^4,\s_4,\psi,\xi)$.
\end{rmk}

Recall that there is a group epimorphism $\pi:\s_4\to\s_3$ with kernel
$H=\langle(12)(34),(13)(24),(23)(14)\rangle$. Moreover, $\pi(\mO_2^4)=\mO_2^3$.
Let $\mQ$ be one of the ql-data from Subsection \ref{pointe-over-sn}, for $n=4$.

\begin{lema}\label{lem:ql datum}
Let $\mQ$ as above. Take $\gamma=0$ if $\mQ=\mQ^{-1}_4$. Then there is an
epimorphism of algebras $\hq\twoheadrightarrow\mH(\mQ_3^{-1}[\lambda])$.
\end{lema}
\pf
Consider the ideal $I$ in $\hq$ generated by the element $H_{(12)}H_{(34)}-1$,
and let $\mL=\hq/I$. We have
\begin{align*}
&H_{(14)}H_{(23)}=\ad(H_{(24)})(H_{(12)}H_{(34)}) && \text{ so } &&
H_{(14)}H_{(23)}=1 && \text{in }\mL,\\
&a_{(34)}=\ad(H_{(14)}H_{(23)})(a_{12}) &&  \text{ so } && a_{(34)}=a_{(12)} &&
\text{in }\mL.
\end{align*}
Analogously, $H_{(13)}=H_{(24)}$, $ a_{(14)}=a_{(23)}$ and $a_{(24)}=a_{(13)}$
in $\mL$. Since, for this ql-data, the action $\cdot:\s_4\times X\to X$ is given
by conjugation and $g:X\to \s_4$ is the inclusion, relations
\eqref{eqn:hq1} and \eqref{eqn:hq2} in the definition of $\hq$ are satisfied in
the quotient. It is now easy to check that the quadratic relations
\eqref{eqn:hq3} defining $\hq$ become in the quotient the corresponding ones
defining the algebra $\mH(\mQ_3^{-1}[\lambda])$.
\epf

\begin{prop}\label{pro:4nonull}
Assume that $(Y,F,\psi,\xi)$ satisfies
\begin{itemize}
\item[(i)] $\xi_{C_i}=\xi_{C_j}$, $\forall\,i,j\in Y$.
\end{itemize}
If $\mQ\neq \mQ_4^\chi(\lambda)$ assume in addition that
\begin{itemize}
\item[(ii)] if $i,j\in Y$, $i\rhd j=j$ and $(i,j)\in C$ then $\xi_C=2\xi_i$.
\end{itemize}
Then the algebra $\mA(Y,F,\psi,\xi)$ is not null.
\end{prop}
\begin{proof}
Assume first that $\psi\equiv1$. Now, given a datum $(Y,F,\psi,\xi)$,
$\pi(F)<\s_3$ and it is easy to see that $\pi(Y)$ is a subrack of $\mO_2^3$.
Moreover, it follows that $\xi$ is compatible with the triple
$(\pi(Y),\pi(F),\psi)$. Then we have the algebra $\mA(\pi(Y),\pi(F),\psi,\xi)$.
As in Lemma \ref{lem:ql datum}, it is easy to see that if we quotient out by the
ideal generated by $\langle e_fe_g\,:\,fg^{-1}\in N\rangle$, then we have an
algebra epimorphism
$\mA(Y,F,\psi,\xi)\twoheadrightarrow\mA(\pi(Y),\pi(F),\psi,\xi)$. As these
algebras are non-zero by Proposition \ref{pro:3nonull}, so is
$\mA(Y,F,\psi,\xi)$.

Notice that in the case in which $(Y,F,\psi,\xi)$ is associated with the
ql-datum $\mQ_4^\chi(\lambda)$, assumption (ii) is not needed, since the first equation in Definition \ref{def:compatible}
implies that, if $i,j\in Y$ are such that $i\rhd j=i$ and $C\in\mR'$ is the
corresponding class, then $\xi_C=0$ and this relation is contained in the ideal
by which we make the quotient.

The case $\psi\neq1$ follows now as in the proof of Proposition
\ref{pro:3nonull}.
\end{proof}

\end{document}